\documentclass[a4paper,11pt]{article}
\usepackage[latin1]{inputenc}
\usepackage[english]{babel}
\usepackage{amssymb, amsmath,amscd,amsxtra,amsfonts, euscript,latexsym, dsfont, mathrsfs,
eufrak, yfonts, longtable,  textcomp, marvosym}
\usepackage{graphics,graphicx}
\usepackage{epsf}
\usepackage{enumerate}
\usepackage{theorem}
\usepackage{fullpage}

\usepackage{array}

\usepackage{
mathptmx,
newcent,
}

\usepackage{hyperref}
\hypersetup{
    backref=true, 
    pagebackref=true,
    plainpages=false,
    bookmarks=true,         
    unicode=false,          
    pdftoolbar=true,        
    pdfmenubar=true,        
    pdffitwindow=true,      
    pdftitle={My title},    
    pdfauthor={Author},     
    pdfsubject={Subject},   
    pdfcreator={Creator},   
    pdfproducer={Producer}, 
    pdfkeywords={keywords}, 
    pdfnewwindow=true,      
    colorlinks=true,       
    linkcolor=blue,          
    citecolor=blue,        
    filecolor=blue,      
    urlcolor=magenta,          
    urlbordercolor=0 1 1
}

\newtheorem{thm}{Theorem}[section]
\newtheorem{defi}[thm]{Definition}

\newtheorem{prop}[thm]{Proposition}
\newtheorem{ex}[thm]{Example}
\newtheorem{oss}[thm]{Remark}
\newtheorem{coro}[thm]{Corollary}
\newtheorem{conj}[thm]{Conjecture}

\newcommand{\rn}{\mathbb{R}^n}
\newcommand{\dimo}{\noindent{\bf Proof : }}
\newcommand{\Min}{\operatorname{Min}}
\newcommand{\res}{\mathop{\hbox{\vrule height 7pt width .5pt depth 0pt
\vrule height .5pt width 6pt depth 0pt}}\nolimits}
\newcommand{\lt}{\{u > t\}}

\newcommand{\boss}{\begin{oss}\rm }
\newcommand{\eoss}{\end{oss}}
\newcommand{\bex}{\begin{ex}\rm }
\newcommand{\eex}{\end{ex}}

\newcommand{\qed}{\thinspace\null\nobreak\hfill\hbox{\vbox{\kern-.2pt\hrule
height.2pt depth.2pt\kern-.2pt\kern-.2pt \hbox to2.5mm{\kern-.2pt\vrule
width.4pt \kern-.2pt\raise2.5mm\vbox to.2pt{}\lower0pt\vtop to.2pt{}\hfil
\kern-.2pt \vrule width.4pt\kern-.2pt}\kern-.2pt\kern-.2pt\hrule
height.2pt depth.2pt \kern-.2pt}}\par\medbreak}

\newcommand\supspace{\rule{0pt}{3.6ex}}

\def\BV{{\mathrm{BV}}}
\def\H{{\mathbf{H}}}
\def\mdiv{\operatorname{div}}

\def\ds{\displaystyle}
\def\oF{\overline{F}}
\def\oW{\overline{W}}

\def\R{\mathbb{R}}

\def\N{\mathbb{N}}
\def\Om{\Omega}
\def\BV{{\mathrm{BV}}}
\def\SBV{{\mathrm{SBV}}}

\def\wpq#1#2{{{\mathrm W}^{#1,#2}}}
\def\lp#1{{\mathrm L}^{#1}}
\def\mdiv{\operatorname{div}}
\def\Ind{\operatorname{I}}
\def\Int{\operatorname{Int}}
\def\mod{\operatorname{mod}\;}
\def\cont{{\mathrm{C}}}

\def\cont{{\mathrm{C}}}

\def\hausp#1{{\cal H}^{{#1}}}

\def\haushn{\hausp{n-1}}
\def\one#1{\mathds{1}_{#1}}
\def\vari{{\displaystyle{\mathbf{v}}}}

\renewcommand{\sectionmark}[1]%
{\markright{\MakeUppercase{\thesection.\#1}}}

\title{Gradient Young measures, varifolds, and a  generalized Willmore functional}

\author{Simon Masnou\footnote{Universit\'e de Lyon, CNRS UMR 5208, Universit\'e Lyon 1, Institut Camille Jordan, 43 bd du 11 novembre 1918, F-69622 Villeurbanne cedex, France. Email: {\tt masnou@math.univ-lyon1.fr}}\quad and\quad 
Giacomo Nardi\footnote{Universit\'e Pierre et Marie Curie Paris 6, CNRS UMR 7598, Laboratoire Jacques-Louis Lions, F-75005, Paris,
France. Email: {\tt nardi@ann.jussieu.fr}} }

\date{August 31, 2012}

\begin{document}

\maketitle

\begin{abstract}
 Being $\Omega$ an open and bounded Lipschitz domain of $\R^n$, we consider the {\em generalized Willmore functional} defined on $\lp{1}(\Om)$ as
 $$F(u)=\left\{\begin{array}{ll}
\ds\int_{\Om}|\nabla u|(\alpha+\beta|\mdiv
\frac{\nabla u}{|\nabla u|}|^p)\,dx&\text{if $u\in \cont^2(\Om),$}\\
+\infty&\text{else},\end{array}\right.$$
where $p>1$, $\alpha>0$, $\beta\geq 0$. We  propose a new framework, that combines varifolds and Young measures, to study the relaxation of $F$ in $\BV(\Omega)$ with respect to the strong topology of  $\lp{1}$. 
\end{abstract}

\section{Introduction}
Let $\Omega$ be an open bounded Lipschitz domain of $\mathbb{R}^n$. We address in this paper the problem of identifying the relaxation (with respect to the strong topology of $\lp{1}(\Omega)$) of the functional
$$F(\cdot,\Omega) : \;u\in\BV(\Omega) \mapsto
\left\{
\begin{array}[m]{ll}
\displaystyle{\int_{\Omega}\left|\nabla u\right|\left(\alpha+\beta\left|\mdiv\frac{\nabla u}{\left|\nabla u\right|}\right|^p\right)dx} & \mbox{\;if\,} u\in
\mbox{C}^2(\Omega)\\
+\infty & \mbox{\;otherwise}
\end{array} \right.
$$
with $p>1$, $\alpha>0$, $\beta\geq 0$ and the convention that the integrand is 0 wherever $|\nabla u| = 0$.  Here, $\BV(\Omega)$ denotes the space of functions of bounded variation in $\Omega$, see~\cite{AFP}. Without loss of generality and to simplify the notations, we shall assume in the sequel that $\alpha=\beta=1$.

This functional appears, under various forms, in the context of optimal design of shapes or digital surfaces in 3D~\cite{Ballester}, modeling and approximation of elastic membranes, or folding in multi-layered materials~\cite{bath1691}, image or surface processing~\cite{MasnouMorel,MM,ChanKangShen,Ballester}. In particular, it has been introduced in~\cite{MasnouMorel,MM} as a variational model in the context of digital image inpainting, i.e. the problem of recovering an image that is known only out of a given domain.  It is also related to a model of amodal completion in a neurogeometric description of the visual cortex~\cite{CittiSarti}.

The functional $F$ has a strong geometric meaning. Indeed, by the coarea formula~\cite{Evans-Gariepy-92,AFP},
\begin{equation}\label{coarea}
F(u, \Omega)=\int_{\mathbb{R}}\left[\int_{\partial \{u>t\}\cap\Omega}\left(1+\left|\mathbf{H}_{\partial\{u>t\}\cap \Omega}\right|^p\right)\;{d} \mathcal{H}^{n-1}\right]\;{d} t  \quad\quad \forall u\in \cont^2(\Omega) 
\end{equation}
where, for a.e. $t$, $\mathbf{H}_{\partial\{u>t\}\cap\Omega}(x)=-(\mdiv\frac{\nabla u}{\left|\nabla u\right|})\frac{\nabla u}{\left|\nabla u\right|}(x)$ is the mean curvature vector at a point $x\in\partial\{u>t\}\cap\Omega$, and $\hausp{n-1}$ is the $(n-1)$-Hausdorff measure. We call $F$ a {\it generalized Willmore functional} for it naturally relates to the celebrated Willmore energy of an immersed compact oriented surface $f:\Sigma\to\R^N$ without boundary, defined as
$${\cal W}(f)=\int_\Sigma|{\mathbf H}|^2dA$$
 with $dA$ the induced
area metric on $\Sigma$. 

\par Minimizing $F$ (for instance under fat boundary constraints) raises immediate difficulties for a simple reason: the functional is not lower semicontinuous with respect to the strong convergence in $\lp{1}$, as can be seen immediately from the following classical example~\cite{BDP}:
\begin{ex}\label{ex:1}{\rm
Being $E$ and $\Omega$ the planar sets drawn on Figure~\ref{cusp12}, left, let $u=\mathds{1}_E$. Obviously, $u\in \BV(\Omega)$ can be approximated in $\lp{1}$ by a sequence $\{u_h\}\subset \cont^2_{\mathrm c}(\Omega)$ of functions with isolevel lines as in Figure~\ref{cusp12}, right.
\begin{figure}[h]
\begin{center}
\includegraphics[height=3cm]{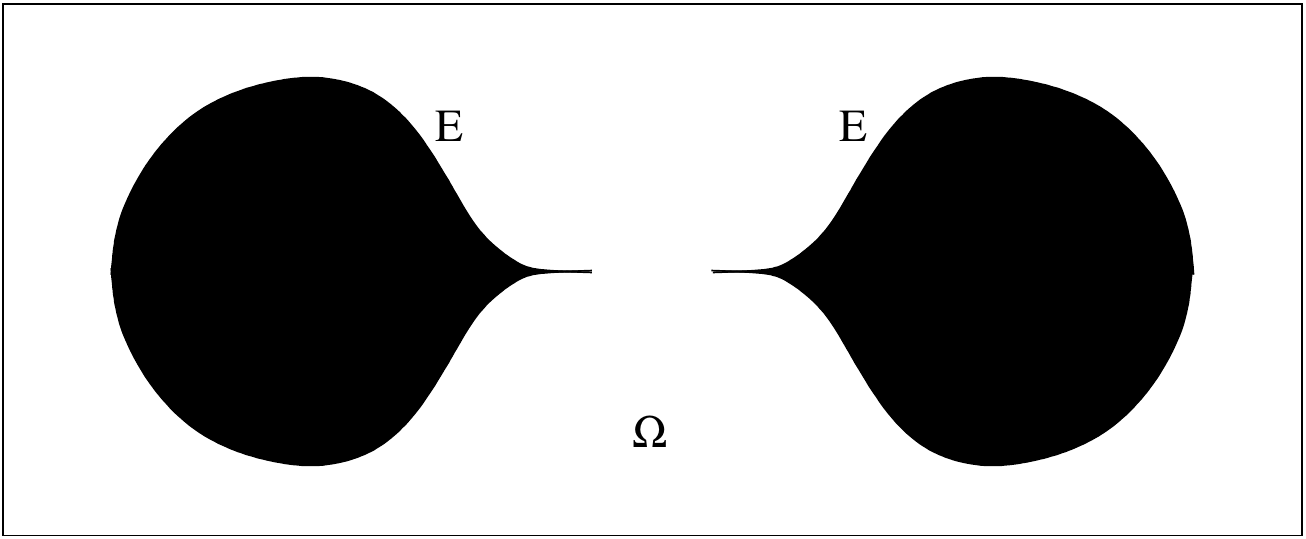}\includegraphics[height=3cm]{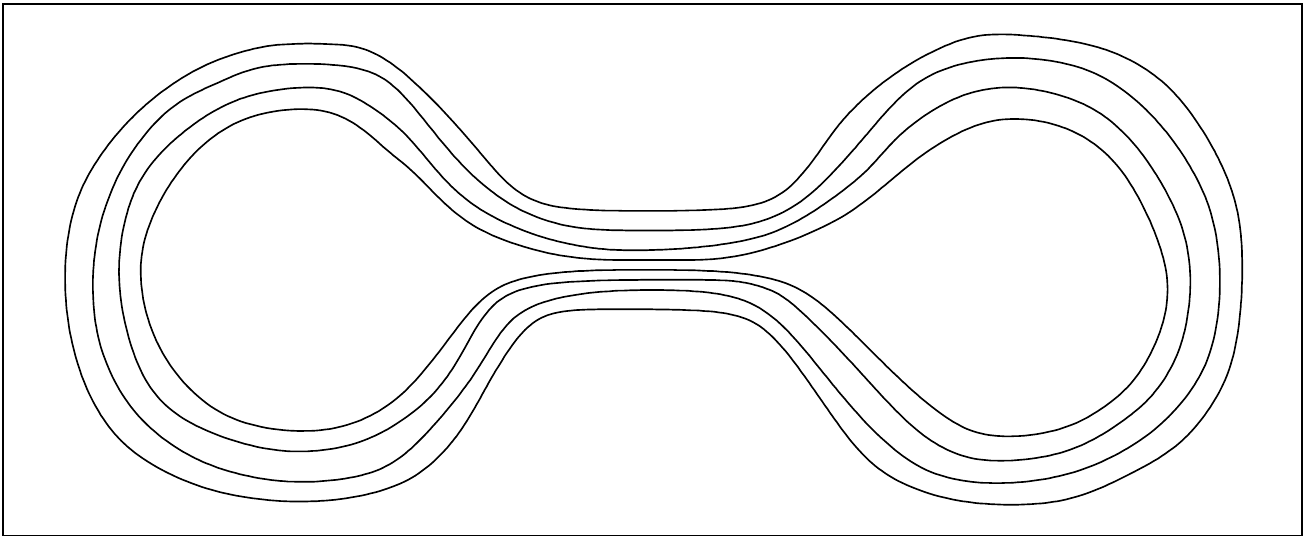}
\caption{Left: $u=\mathds{1}_E$ with $\overline{F}(u,\Omega)<\infty$. Right: Isolevel lines of a smooth approximating function.}\label{cusp12}
\end{center}
\end{figure}
It is easy to check that $\underset{h \rightarrow \infty}{\liminf}\; F(u_h) <\infty$ but, since $u\notin \cont^2(\Omega)$, we have $F(u, \Omega)=\infty$ so $F$ is not lower semicontinuous.
}
\end{ex}

The usual technique in calculus of variations to overcome this difficulty consists in relaxing $F$, i.e., introducing the functional
$$ \oF(u,\Omega) = \inf\left\lbrace \underset{h\rightarrow\infty}{\liminf}\;F(u_h, \Omega) : u_h\overset{\lp{1}(\Omega)}{\longrightarrow}u\right\rbrace .$$
As a relaxation, this functional has the interesting property of being lower semicontinuous in $\lp{1}$~\cite{Dal}. Together with the relative compactness of $\BV$ in $\lp{1}$, it guarantees that the infimum of $F$ coincides with a minimum of $\oF$, which somewhat solves the minimization problem. It remains however that not much can be said neither about the minimizers of $\oF$ nor, more generally, about $\oF(u)$ for a general function $u$ with bounded variation.

Partial results have been obtained in~\cite{AM,LM} in the case where $u$ is smooth. Combining the techniques used in these papers with the more recent~\cite{Menne}, it can be proved that, in any space dimension $n$ and for any $p\geq 1$, $F(u)=\oF(u)$ when $u$ is $\cont^2$. What about more general functions?

Examining again the previous example, it is clear that $\oF(u)<+\infty$ since $(u_h)$ has uniformly bounded energy and converges to $u$. Besides, it is equivalent to study $\oF$ for the function $u=\one{E}$ or to study the relaxation at $E$ of the following functional that acts on measurable sets (in our example $n=2$):
$$A\subset\R^n\;\mapsto\;W(A)=\left\{\begin{array}{ll}
\ds\int_{\partial A}(1+ |{\H}_{\partial A}|^p)\;{d} \haushn& \text{if $\partial A$ is smooth}\\
+\infty&\text{otherwise}\end{array}\right.$$
The relaxed functional associated with $W$ is
$$\oW(A)=\inf\{\liminf_{h\to\infty}\,W(A_h),\;(\partial A_h)\;\text{smooth},\;|A_h\Delta A|\to 0\},$$
where $|\cdot|$ denotes the Lebesgue measure. 

\begin{figure}[h]
\begin{center}
\includegraphics[height=2.5cm]{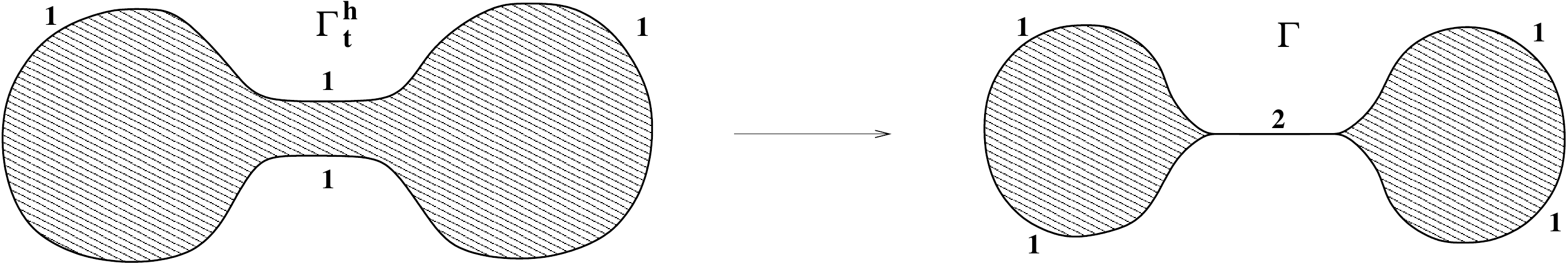}
\caption{Accumulation at the limit of the boundaries of sets that approximate in measure the set $E$ of Figure~\ref{cusp12}.}\label{cusp3}
\end{center}
\end{figure}

The properties of bounded sets $A\subset\R^2$ such that $\oW(A)<+\infty$ and the explicit representation of $\oW(A)$ have been carefully studied in~\cite{BDP,BM1,BM}. Such sets have finite perimeter (by definition of the energy) and, by explicit representation, we mean that $\oW(A)$ can be written in terms of the $\wpq{2}{p}$ norms of  a collection of curves that cover the essential boundary $\partial^*A$ of $A$. This can again be easily understood from Example~\ref{ex:1} and Figures~\ref{cusp12},~\ref{cusp3}: a "good" way to approximate $E$ in measure consists in choosing a set $E_h$ whose boundary $\Gamma^h$ is represented in Figure~\ref{cusp3}. These sets have uniformly bounded energy and, as shown in~\cite{BDP}, $\oW(E)$ coincides with $W(\Gamma)=\int_\Gamma(1+|\kappa_\Gamma|^p)d\hausp{1}$ (with $\kappa_\Gamma$ the curvature along $\Gamma$) where $\Gamma$ is the limit curve represented in Figure~\ref{cusp3}, right, with its multiplicity.

Having in mind the expression \eqref{coarea} of $F$ through the coarea formula, it is natural to expect that, at least in dimension 2, the relaxed energy $\oF(u)$ of a function $u$ of bounded variation can be written in terms of the energies $W(\Gamma_t)$ of systems of curves $(\Gamma_t)$ that cover the essential boundaries $\partial^*\{u>t\}$. This is exactly what happens, as we proved in the companion paper~\cite{MN}, among other results, by generalizing techniques that were proposed in~\cite{MM} for a more restrictive boundary value problem involving the same energy. The precise statements will be recalled in Section~\ref{primoarticolo}. 

The techniques developed in~\cite{BDP,BM1,BM,MM,MN} depend strongly on parameterizations of curves and can hardly be generalized to higher space dimensions. Indeed, in dimension strictly greater than 2, parameterizations of hypersurfaces are much harder to handle in our context especially since the energy of interest controls the mean curvature vector only. We will come back later on this issue, that was the main motivation for the new framework that we introduce in this paper and that involves two specific tools: Young measures, that play a fundamental role in many problems of the calculus of variations, and varifolds, that appear to be very useful to handle generalized surfaces and a weak notion of mean curvature.

\paragraph{Varifolds}
The basic idea behind rectifiable varifolds, that will be introduced with more details in Section~\ref{teoriavarifold}, is that each rectifiable $k$-subset $M\subset\R^n$ can be endowed with a multiplicity function $\theta_M$ and associated with the measure $\theta_M\hausp{k}\res M$. 
The associated varifold is the Radon measure $V_M=\vari(M,\theta_M)=\theta_M\hausp{k}\res M\otimes\delta_{T_M(x)}$ on the product space $G_k(\Omega)=\Omega\times G(n,k)$, with $G(n,k)$ the Grassmaniann of $k$-subspaces in $\R^n$ and $T_M(x)$ the tangent space to $M$ at $x\in\Gamma$. Therefore, varifolds carry information both on spatial localization and tangentia behavior. Varifolds have nice properties, among which the possibility to use a weak notion of mean curvature, the continuity of the mass, a compactness property, the lower semicontinuity of some useful second-order energies, etc. They are actually much more adapted than parameterizations for handling sequences of $k$-surfaces in $\R^n$ when there is a control on the mean curvature.

Denoting $\mu_V(\cdot)=V(\cdot\times G(n,k))$ the mass of a $k$-varifold $V$ in $\R^n$, we show in the following table how notions that are naturally defined for smooth $k$-sets can be easily translated in terms of $k$-varifolds. Here, $X$ denotes a smooth vector field with compact support.
\begin{center}
\begin{tabular}{|m{3cm}|m{7cm}|m{4.5cm}|}
\hline
                & $M$ closed, smooth $k$-set                                &  $V$ $k$-varifold                                       
                \vspace*{0.2cm}\\
                \hline 
\supspace Mass            & $\mathcal{H}^{k}(M)$                                 &  $\mu_V=V(\Omega\times G(n,k))$                      \vspace*{0.2cm}\\
\hline
\supspace First variation & $\ds\int_{M} {\rm{div}}_{M} X d\hausp{k}$  &  \supspace $\delta V(X)\!\!=\!\!\ds\int_{G_{k}(\Omega)}\!\!\!\!\!\! {\rm{div}}_S X dV(x, S)$ \vspace*{0.2cm}\\
\hline
\supspace Mean curvature   vector    & ${\bf H}_M$ s.t $\ds\int_{M} {\rm{div}}_{M} X d\hausp{k}\!\!=\!\!-\!\!\int_M {\bf H}_M\cdot Xd\hausp{k}$                             &  $\ds {\bf H}_V=-\frac{\delta V}{\mu_V}$\\
\hline
\end{tabular}  
\end{center}

In particular, this table shows the divergence theorem that relates (the integral of) the tangential divergence of a smooth vector field with (the integral of) the mean curvature vector. Recalling that, for a smooth function $u$,  $-(\mdiv\frac{\nabla u}{|\nabla u|})\frac{\nabla u}{|\nabla u|}(x)$ is the mean curvature vector at $x\in \partial\{u>t\}$, and
$$\mdiv_{\{y,\,u(x)=u(x)\}}X=\mdiv_{{\nabla u^\perp}}X$$
we calculate
$$\begin{array}{lll}
\ds\int_{\Omega}|\nabla u|\mdiv_{{\nabla u^\perp}}X\,dx&=&\ds\int_{-\infty}^{+\infty}\int_{\Omega\cap \{u=t\}}\mdiv_{{\nabla u^\perp}}X\,dx\,dt=-\ds\int_{-\infty}^{+\infty}\int_{\Omega\cap \{u=t\}}
{\bf H}_{ \{u=t\}}\cdot X\,dx\,dt\\[5mm]
&=&-\ds\int_{\Omega}|\nabla u|{\bf H}_{ \{y,\,u(y)=u(x)\}}\cdot X\,dx\end{array}$$
This looks exactly like the formula provided by the divergence theorem, except that the Hausdorff measure has been replaced by the measure $|\nabla u|dx$. This observation is the core of our approach in a less regular context: roughly speaking, given a function $u$ of bounded variation, we will define a varifold associated with the mass represented by the total variation. Then the first variation of the mass can be computed (like above), and considering the measure provided by the Riesz representation theorem, its Radon-Nikodym derivative with respect to the mass finally yields the mean curvature.

The first delicate issue is to extend properly to $\BV$ the quantity $\ds\int_{\Omega}|\nabla u|\mdiv_{{\nabla u^\perp}}X\,dx$ that belongs to the general family of mappings $\ds u\mapsto \int_\Omega f(x,\nabla u)dx$. Studying such mappings in $\BV$ is the purpose of~\cite{KR}, where suitable tools are defined based on the theory of Young measures.

\paragraph{Young measures}
They were introduced by L.C. Young~\cite{Y1,Young2,Young3} to describe limits of minimizing sequences for integrals of the type
$$\int f(x,u)\,dx\quad \mbox{or}\quad \int f(x,u,\nabla u)\,dx$$
Young measures are particularly useful when classical minimizers do not exist. They can  handle complex situations with concentration, oscillation, or diffusion phenomena. They find many applications in calculus of variations, optimal control theory, optimal design, variational modeling of nonlocal interactions, etc.~\cite{VAL,CFV,Ped}.

The typical situation where they arise is the following: $\Omega$ being bounded, take a sequence $(v_h)$ that converges weakly to $v$ in $\lp{\infty}(\Omega,\R^n)$, and look at $f(x,v_h(x))$ with $f$ continuous and nonlinear. Then a classical theorem due to L.C. Young states that there exists a family of probability measures $(\nu_x)_{x\in \Omega}$, called Young measure generated by the sequence $(v_h)$,  such that 
$$\int_{\rn}z \;{d} \nu_x(z) = v(x) \quad \quad \mathcal{L}^n-\text{a.e. }x$$
and, up to a subsequence,
$$ \int_{\Omega} f(x, v_h) \;{{d} x}=\int_\Om\int_{\R^n}f(x,z)\delta_{\nu_h}(z)\,dx \rightarrow \int_{\Omega}\int_{\rn} f(x, z) \nu_x(z)\;{{d} x}.$$
In other words, the impossibility to use the continuity of $f$ is overcome by introducing a measure that, in some sense, carries the information out of $f$. 

\par A frequent situation in the calculus of variations concerns the case where $v_h$ are gradients, i.e. $v_h = \nabla u_h$ and $v=\nabla u$ for some $u_h, u \in\wpq{1}{p}(\Omega)$. As above, every sequence of gradients that weakly converges in $\lp{p}$ generates a Young measure, called gradient Young measure. 
\par  Then it is natural to ask which  families of probability measures are generated by  sequences of gradients or, in other words, can one characterize the set of gradient Young measures?
\par   In \cite{KP1,KP2} the authors study the gradient Young measures generated by a sequence of gradients converging weakly in $\lp{p}(\Omega, \mathbb{R}^m)$ ($p>1$)  and their  characterization essentially depends on the condition 
$$\int_{\Omega}\int_{\rn}|z|^p \;{d} \nu_x(z){d} x < \infty. $$
Their results are generalizable to $p=\infty$ and to the vectorial case (i.e. for $\mathbb{R}^d$-valued functions, with $d>1$) , see~\cite{Ped} for precise statements.
\par In the applications, if $p>1$, the weak convergence follows from a uniform bound on the $\wpq{1}{p}$ norm of the gradients, but in the case $p=1$ the space $\wpq{1}{1}$ is not reflexive so, to infer weak relative compactness, the sequence $\{\nabla u_h\}$ should be equi-integrable, which is hard to establish in the applications. As an alternative, the weak-* topology of $\BV(\Omega)$ can be considered, and leads to an extension of the concept of Young measures.
\par In \cite{PM,AB97,KR}, a new formulation for Young measures is introduced to extend the classical theory to the framework of functions of bounded variation. A generalized gradient Young measure $\nu$ is defined as a triplet of measures $\nu=(\nu_x, \nu_x^{\infty},\lambda_\nu)$ where $(\nu_x)_{x\in \rn}$ is a family of probability measures on $\rn$, $\lambda_\nu$ is a positive bounded Radon  measure on $\overline{\Omega}$ and $(\nu_x^{\infty})_{x\in \rn}$ is a family of probability measures on $\mathbb{S}^{n-1}$, the unit sphere of $\rn$.
The Young measure representation is extended to
$$\langle\langle\nu,f\rangle\rangle:=
\int_{\Omega}\int_{\rn}f(x,z)\;{d} \nu_x(z)\;{{d} x} +
\int_{\overline{\Omega}}\int_{\mathbb{S}^{n-1}}f^{\infty}(x,z)\;{d} \nu_x^{\infty}(z)\;{d} \lambda_{\nu}(x)
$$
where $f^{\infty}$ is the recession function defined as $\ds f^{\infty}(x,z) := \underset{\underset{\underset{t\rightarrow
\infty}{z'\rightarrow z}}{x'\rightarrow x}}{\lim}\dfrac{f(x', tz')}{t}.$
\par In \cite{KR} a characterization theorem for generalized Young measures is proved. As in the case of classical Young measures, a necessary condition in the real-valued case for having a generalized gradient Young measure (i.e. generated by a sequence of gradients of functions that converge  weakly-* in $\BV$) is
$$\int_{\Omega}\int_{\rn}|z| \;{d} \nu_x(z){d} x  + \lambda_\nu(\overline{\Omega})< \infty $$
 We refer to \cite{KR} for general results in the vectorial context of $\BV(\R^n,\R^m)$.  In Section \ref{Youngmeasure} we recall the  main results in the real-valued case $m=1$ .
\par The examples given in \cite{AB97,KR}, and those from Section \ref{Youngmeasure}, show that Young measures are totally suitable for describing the concentration and oscillations effects generated by the weak convergence of  gradients. In fact, limit Young measures contain analytic and geometric information; they depend on the converging sequence (and not only on its weak limit!) so they carry some information about the weak limit of the sequence of  gradients and the intrinsic features of the sequence.
\paragraph{Young varifolds}
We have now the material to introduce the Young varifolds, i.e. a suitable class of varifolds generated by Young measures which allows us to formalize our problem in the varifolds framework.~\footnote{M. Novaga kindly brought to our attention, while the current paper was in the final correction phase, the reference~\cite{lorentz} where a generalization of Almgren's theory of varifolds in a Lorentzian setting is proposed. In a different context and for different purposes, it shares with our work the idea of disintegrating and indexing the measures that we borrowed from~\cite{KR} while it is done using ad-hoc varifolds in~\cite{lorentz}.}
\par For every $f\in \cont_{\mathrm c}(G_{n-1}(\Omega))$ let
$$g:\;(x,z)\in \Omega \times \rn \mapsto g(x,z)= |z|f(x, z^\perp)$$
where $z^\perp$ is the element of $G_{n-1}(\Omega)$ perpendicular to $z$. It is easy to check that for every $k\in \mathbb{R}$ we have $(kz)^\perp = z^\perp$ (as elements of $G(n, n-1)$) so we get $g^\infty(x,z)= f(x,z^\perp).$
\par A varifold $V$ is a Young varifold if there exists a Young measure $\nu$ such that
\begin{equation}\label{def:youngvari}
\int_{G_{n-1}(\Omega)} f(x,S) \;{d} V_\nu(x,S) = \int_{\Omega}\int_{\rn}|z|f(x,z^\perp) \;{d} \nu_x\;{{d} x} + \int_{\Omega}\int_{\mathbb{S}^{n-1}}f(x,z^\perp) \;{d} \nu_x^\infty\;{d} \lambda_\nu
\end{equation}
for every $f\in \cont_{\mathrm c}(G_{n-1}(\Omega))$. $V=V_\nu$ is called the Young varifold associated to $\nu$. 
\par The definition of a Young varifold is particularly explicit for smooth functions. If $u\in\cont^2(\Omega)$ we consider the Young measure $\nu=(\nu_x,\nu_x^\infty,\lambda_\nu)$ defined by
$$\nu_x = \delta_{\nabla u_x}\;, \quad \nu_x^\infty=0\;,\quad \lambda_\nu=0$$
and it follows that
$$ \int_{G_{n-1}} f(x,S) \;{d} V_\nu(x,S) = \int_{\Omega} |\nabla u|f(x, (\nabla u)^\perp(x))\;{{d} x}.$$
The mass of the varifold is defined by $$\mu_{V_{\nu}}(E)=V_{\nu}(G_{n-1}(E))=\int_{E}\int_{\rn}|z|
\;d \nu_x(z)\;{d x}+\lambda_{\nu}(\overline{E}) \quad  \; \forall E\subseteq \Omega.$$
and the first variation is
$$\delta V_{\nu}(X)=
\ds \int_{G_{n-1}(\Omega)}\mdiv_S X(x)\; \;d V_{\nu}(x,S) =\begin{array}[t]{l}\ds\int_{\Omega}\int_{\rn}|z|\mdiv_{ z^\perp}X\;\;d \nu_x(z)\;{d x}+ \ds\int_{\Omega}\int_{\mathbb{S}^{n-1}}\mdiv_{z^{\perp}}X\;\;d \nu_x^{\infty}(z)\;d \lambda_{\nu}\end{array}$$

We can now, as we did for sets, show how the usual notions for smooth functions can be extended to the framework of Young varifolds.

\begin{center}
\begin{tabular}{|m{3cm}|m{4cm}|m{6cm}|}
\hline
                & $u\in\cont^{2}$                                              &  $V_\nu$ Young varifold                                     \vspace*{0.2cm}\\
\hline\vspace*{0.4cm}
Mass            & $ |D u|$                                                            &  $\mu_{V_\nu}$                      \\
\hline
 First variation & $\ds\int_{\mathbb{R}^2} |\nabla u|{\rm{div}}_{\nabla u^\perp} X dx$      & \supspace $\delta V_\nu(X)=\ds\int_{G_{n-1}(\Omega)} {\rm{div}}_S X dV_\nu(x, S)$ \\
\hline
\supspace Mean curvature vector    & $-\ds( {\rm{div}}\frac{\nabla u}{|\nabla u|}) \frac{\nabla u}{|\nabla u|}$    &  $\ds {\bf H}_V=-\frac{\delta V_\nu}{\mu_{V_\nu}}$\\
\hline
\end{tabular}
\end{center}
Finally, we define the generalized Willmore energy associated with a Young varifold as
$$W(V)=\int_{\Omega}\left( 1+ \left|{\bf H}_V\right|^p \right)
\;d \mu_{V}.$$

\par The paper is devoted to defining carefully Young varifolds, exhibiting some of their properties and investigating the relationship between Young varifolds and the relaxation problems for $F$ and $W$. Given a function $u\in\BV$ such that $\oF(u)<\infty$, we focus on the class $\mathbb{V}(u)$ of all Young varifolds $V_\nu$ such that $\|\delta V_\nu\| <\!< \mu_{V_\nu}$ and $\nu$ is the limit of gradient Young measures $\nu_{\nabla u_n}$ where $(u_n)$ converges  weakly-* to $u$ in $\BV$. Studying the Young varifolds in  $\mathbb{V}(u)$ is somewhat delicate. In {\it dimension ${\mathit 2}$}, using the results of~\cite{MN}, we prove (Theorems~\ref{sistemaYoung} and~\ref{sistemaYoung}, Corollary~\ref{phiyoung}) that for every $u\in\SBV$ with compact support and such that $\oF(u)<\infty$, there exists $V\in \mathbb{V}(u)$ such that $W(V)=\oF(u)$. We conjecture that a similar result holds in higher dimension but the proof remains so far out of reach. We shall comment this point later on. So far, we are able to prove in Theorem~\ref{minVu} that, {\it in any dimension${\mathit \geq 2}$},
\begin{equation}
\overline{F}(u,\Omega)\geq \Min\; \{W(V): V\in \mathbb{V}(u)\},\label{eq:ineq}
\end{equation}
\par There is no hope that equality holds in general in~\eqref{eq:ineq} as arises from simple two-dimensional examples, see Remark~\ref{corocomp} and Proposition~\ref{propcex}. Therefore, what additional assumptions must be taken in $\mathbb{V}(u)$ to guarantee the equality? Clearly, in dimension $2$, it follows from~\ref{sistemaYoung} that a necessary assumption is the existence of a tangent {\it everywhere} on the support of the concentration measure $\lambda$, see also~\cite{BM}. In contrast, things are really unclear in higher dimension and are the purpose of ongoing research. It follows from the results of Menne in~\cite{Menne} that $\lambda$ can be decomposed into $(n-1)$--fibers whose supports are $\cont^{2}$-rectifiable, but this regularity remains too weak even in dimension $2$ (a tri-segment is $\cont^{2}$-rectifiable but there is no tangent at the triple point). A more accurate characterization is needed, which has of course to do with the largely open problem of characterizing the boundaries of $n$-sets, $n\geq 3$, whose relaxed Willmore energy is finite.

\medskip
\par Why do we believe that Young varifolds are the right tools for tackling the problem of representing $\oF$ in dimension higher than $2$? Because they offer the possibility to carry all together and implicitly the concentration at the limit of the boundaries $\partial\{u\geq t\}$, using a unique representation of the general form $(\nu_x,\nu_x^\infty,\lambda)$, and because compactness and semicontinuity of the energy under constraints are obtained very easily (see Corollary~\ref{corocomp}). As for the information carried by Young varifolds, it must be emphasized from~\eqref{def:youngvari} that a Young varifold of $\mathbb{V}(u)$ for $u\in\BV(\R^n)$ is a $(n-1)$--varifold (it acts on $G_{n-1}(\Omega)$) but is not necessarily rectifiable: the support in $\R^2$ of the Young varifold associated with $u(x,y)=x$ is the whole plane. A Young varifold in $\mathbb{V}(u)$  must rather be seen as a {\it fiber bundle} whose fibers are $(n-1)$-rectifiable varifolds.

\medskip
\par Are there alternative approaches to the problem? It has been shown in~\cite{AM} that the study of $\oF(u)$ when $u$ is {\it smooth} can be tackled considering explicitly all boundaries $\partial\{u\geq t\}$. Following~\cite{AM,MN}, let us examine whether the same strategy is applicable when $u$ is possibly unsmooth. Take $u\in\lp{1}(\R^n)$ and a sequence of smooth functions $(u_h)$ converging to $u$ in $\lp{1}(\R^n)$ and such that $F(u_h)\to \oF(u)$ as $h\to\infty$. Possibly extracting a subsequence, one can assume that for almost every $t$, $\{u_h>t\}$ converges to $\{u>t\}$ in measure. In addition, by Fatou's
Lemma,
$$\int_\R\liminf_{h\to\infty}\,W(\{u_h>t\})dt\leq\liminf_{h\to\infty}\int_\R W(\{u_h>t\})dt=\oF(u)$$
therefore $\liminf_{h\to\infty}W(\{u_h>t\})<\infty$ is finite for almost every $t$. It follows that, for almost every $t$, the $(n-1)$-dimensional varifolds with unit multiplicity $\vari(\partial\{u_h>t\},1)$ form a sequence with uniformly bounded mass, and uniformly bounded curvature in $\lp{p}$. Since $p>1$, by the properties of varifolds~\cite{Si} and the stability of absolute continuity (see Example 2.36 in~\cite{AFP}), there exists a subsequence $\vari(\partial\{u_{h_k}>t\},1)$ \textit{depending on $t$} and a limit integral $(n-1)$-varifold $V_t$ such that 
$$\int_{\R^n}(1+|{\bf H}_{V_t}|^p)d\|V_t\|\leq \liminf_{h\to\infty}\,W(\{u_h>t\})$$
In addition, one can prove~\cite{AM} that the support $M_t$ of $V_t$ contains $\partial^*\{u>t\}$ for almost every~$t$. What could be the remaining steps to get a representation of $\oF(u)$?
\begin{enumerate}
\item show, if possible, that the limit varifolds $V_t$ are nested, i.e. $\operatorname{int}V_t\subset\operatorname{int}V_{t'}$ if $t>t'$, where $\operatorname{int}V_t$ denotes the set enclosed (in the measure-theoretic sense) by the support of $V_t$. Again, observe that $\ds\int_{\R^n}(1+|{\bf H}_{V_t}|^p)d\|V_t\|\leq \liminf_{h\to\infty}W(\{u_h>t\})$.
\item build a sequence of smooth sets $E_h^t$ (for a suitable dense set of values $t$) such that $\partial E_h^t\to M_t$ (being $M_t$ the support of $V_t$) and $W(E_h^t)\to\ds\int_{\R^2}(1+|{\bf H}_{V_t}|^p)d\|V_t\|$. The varifolds $V_t$ being nested, one could actually build $E_h^t$ so that $E_h^t\subset E_h^{t'}$ if $t>t'$.
\item by a suitable smoothing of the sets $E_h^t$, build a smooth function $\tilde u_h$ such that $F(\tilde u_h)\leq \ds\int_\R W(E_h^t)dt+\frac 1 h$. 
\item passing to the limit, possibly using a subsequence, show that $\tilde u_h$ tends to $u$ in $\lp{1}$ and using the lower semicontinuity of $\oF$, conclude that
$$\oF(u)=\int_\R\int_{\R^n}(1+|{\bf H}_{V_t}|^p)d\|V_t\|\,dt$$
\end{enumerate}
which would be a nice representation formula. The delicate steps in this tentative proof are steps 1 and 2. It is in particular not clear at all whether the limit varifolds are nested. It would be an easy consequence of the existence of a subsequence $(u_{h_k})$ such that the varifolds $\vari(\partial\{u_{h_k}>t\},1)$ converge to $V_t$ for almost every $t$. But this is \textit{false} in general as shown by a counterexample communicated to us by G. Savar{\'e} and fully described (Example 1.2) in the companion paper~\cite{MN}. The example shows a sequence $\{\tilde{u}_n\}\subset \cont^{0}([0,1]^2)$ with smooth level lines $\{\tilde u_n=t\}$ satisfying 
$$\sup_n \int_{\R}  \int_{\partial \{\tilde{u}_n(x)>t\}\cap (0,1)^2}(1+\left|\kappa_{\partial \{\tilde{u}_n(x)>t\}\cap (0,1)^2}\right|^p ) \,d\mathcal{H}^{1} \,dt <  \infty,$$
but such that there exists no subsequence $(t\mapsto\vari(\partial\{u_{h_k}>t\},1))_k$ converging for almost every $t$ to a limit varifold $V_t$.

\par In the particular case of dimension $2$, we overcame this subsequence issue in~\cite{MN} using the fact that the varifolds are supported on $\wpq{2}{p}$ parametric curves. Then step 1 follows from the selection of countably many ``shepherd'' curves, that guide the remaining others, and a diagonal extraction argument that uses the $\cont^{1}$ convergence of the parametric curves. Having parametric curves is crucial for the smoothing step (step 2), and more precisely for moving apart the curves while controlling the energy.
\par Is the same strategy applicable to dimension greater than $2$?  The martingale argument that we used in~\cite{MN} for the diagonal extraction is valid in any dimension, and the convergence of countably many $\cont^{1}$ curves can be replaced by the convergence of countably many integral varifolds, which, even being much weaker, is enough to obtain the limit structure. However, we do not know so far whether step 2 could be generalized to higher dimension. It follows from Menne's results~\cite{Menne} that the limit varifolds are supported on $\cont^{2}$-rectifiable sets but it is far from being clear how these countable coverings can be smoothly deformed while controlling the energy of the underlying set. Above all, we feel that the understanding of the problem could benefit from using a framework that is lighter than explicit unions of integral varifolds, that provides easily relative compactness and semicontinuity of the energy, and this motivated the introduction of Young varifolds.

\medskip
\par The plan of the paper is as follows: the first sections are dedicated to a careful introduction of all notions that we have roughly described so far. More precisely, in Section \ref{primoarticolo} we recall the main definitions and results obtained in \cite{MN}. Section \ref{funzionalimisure} and \ref{teoriavarifold} are devoted, respectively, to a general class of functionals depending on measures and to the varifold theory. In Section \ref{Youngmeasure} we recall a few facts about Young measures, following \cite{KR}. In Sections \ref{Youngvarifold} and  \ref{Youngvarifold2} we define the Young varifolds and their Willmore energy, and we provide several examples  showing that Young varifolds allow to get information about geometric phenomena, like oscillations and concentration, for minimizing sequences. Lastly, we study in section \ref{rappresentazione} the relationship between $\oF$ and the Willmore functional for Young varifolds.

\begin{center}
 {\bf General notations}
\end{center}

$\rn$ is equipped with the Euclidean norm and we will denote by either $\mathcal{L} ^n$ or $|\cdot|$ the Lebesgue measure on $\R^n$. $\mathcal{H}^k$ is the $k$-dimensional Hausdorff measure. The restriction of a measure $\mu$ to a set $A$ is denoted by $\mu\res A$ and  $\text{spt}\; \mu$ is the support of $\mu$. 
\par For two open sets $E, F\subset \rn$ the notation $E\subset\subset F$ means that $\overline{E}\subset F$ and  $\overline{E}$ is compact.
\par If $X$ is a locally compact separable metric space we denote by $\mathcal{M}(X, \rn)$ the space of $\rn$-valued bounded Radon measures and by $\mathcal{M}^+(X)$, $\mathcal{M}^1(X)$ the spaces of  positive Radon measures and  probability measures, respectively. Moreover, given $\mu\in \mathcal{M}(X, \rn)$ and $ \nu \in \mathcal{M}^+(X)$ we denote by $\frac{\mu}{\nu}$ the derivative of $\mu$ with respect to $\nu$ and the Radon-Nikodym decomposition of $\mu$ with respect to $\nu$ is  $\mu= \mu^a +\mu^s= \frac{\mu}{\nu}\nu + \mu^s$.
\par $\cont_{\mathrm c}$, $\cont^{k}, \lp{p}, \wpq{k}{p}, \BV, \SBV$ are the usual function spaces. For a detailed study of the spaces $\BV$ and $\SBV$ of functions with bounded variation, the reader may refer to~\cite{AFP}. If $\Omega \subset \rn$, we say that $\partial \Omega \in \cont^k$ (resp. $\wpq{k}{p}$) if  we can represent locally its boundary as a graph of a $\cont^{k}$ (resp. $\wpq{k}{p}$) function. In particular, $\Omega$ is called a Lipschitz domain if  $\partial\Omega \in \cont^{0,1}$.

\section{Relaxation by a coarea-type formula in dimension $2$}\label{primoarticolo}

We recall in this section the main results proved in \cite{MN} which will be used in the following. Let us start with the notion of system of curves of class $W^{2,p}$:

\begin{defi}\label{sistemi}
 By a system of curves of class $\wpq{2}{p}$ we mean a finite family $\Gamma=\{\gamma_1,...,\gamma_N\}$ of closed curves of class $\wpq{2}{p}$ (and so $\cont^1$) admitting a parameterization  (still denoted by  $\gamma_i$) $\gamma_i\in \wpq{2}{p}\left([0,1],\R^2 \right)$  with unit velocity.
Moreover, every curve of  $\Gamma$ can have  tangential self-contacts but without crossing and two curves of $\Gamma$ can have  tangential contacts but without crossing. In particular, $\gamma_i'(t_1)$ and $\gamma_j'(t_2)$ are parallel whenever $\gamma_i(t_1)=\gamma_j(t_2)$ for some $i,j\in\{1,...,N\}$ and $t_1,t_2\in[0,1]$.
\par The trace  $(\Gamma)$ of $\Gamma$ is the union of the traces $(\gamma_i)$, and the interior of the system $\Gamma$ is
$$ \Int(\Gamma) = \{x \in \R^2 \setminus (\Gamma) : \Ind(x, \Gamma) = 1 \; \mod 2\},\qquad \mbox{with }\Ind(x, \Gamma) = \sum_{i=1}^N \Ind(x, \gamma_i).$$
The multiplicity  function $\Gamma$ is $\theta_{\Gamma} :  (\Gamma) \rightarrow \mathbb{N}$,  $\theta(z)=\sharp \{\Gamma^{-1}(z)\},$ where $\sharp$ is the counting measure. \\
If the system of curves is the boundary of a set $E$ with $\partial E\in \cont^2$, we simply denote it as $\partial E$.
\end{defi}

\boss Remark that, by  previous definition, every $|\gamma_i'(t)|$ is constant for every  $t\in[0,1]$ so the arc-length parameter 
is given by $s(t)=t L_i$ where $L_i$ in the length of $\gamma_i$. Denoting by $\tilde{\gamma}_i$ 
the curve parameterized with respect to the arc-length parameter we have  
$$s\in[0, L_i] \,,\,\, \tilde{\gamma}_i(s)= \gamma_i(s/L_i)\,,\,\, \tilde{\gamma}_i''(s)=\dfrac{\gamma_i''(s)}{L_i^2}.$$
Now, the curvature ${\mathbf k}$ as a functions of $s$, satisfies ${\mathbf k} = \tilde{\gamma}_i''(s)$, which implies
$$ \int_{0}^{L_i} \left( 1+|\tilde{\gamma}_{i}''(s)|^p\right)\mbox{d}s=\int_{0}^{L_i} \left( 1+|{\mathbf k}|^p\right)\mbox{d}s= \int_{0}^{1} \left( |\gamma_{i}'(t)|+L_i^{1-2p}|\gamma_{i}''(t)|^p\right)\mbox{d}t.$$
Then, the condition $\gamma_i \in \wpq{2}{p}\left([0,1],\R^2 \right)$ implies that $\tilde{\gamma}_i\in \wpq{2}{p}\left([0,L_i],\R^2 \right)$ 
and, for simplicity, in the sequel we denote  by $\gamma_i$  the curve parameterized with respect to the arc-length parameter.
\eoss

In dimension 2, the Willmore functional for a system $\Gamma$ of curves of class $W^{2,p}$ is
$$W(\Gamma)=\sum_{i=1}^N W(\gamma_i)=\sum_{i=1}^N \int_{(\gamma_i)} \left( 1+|{\bf k}_{\gamma_{i}}|^p\right)\;{d} \mathcal{H}^1.$$

\begin{defi}\label{sistlim}
We say that $\Gamma$ is a limit system of curves of class $\wpq{2}{p}$ if $\Gamma$ is the weak limit of a sequence $(\Gamma_h)$ of boundaries of bounded open sets with $\wpq{2}{p}$ parameterizations.
\end{defi}

The following class of curve-valued functions will be used for covering the level lines of a real function.

\begin{defi}\label{A}
Let $\mathscr{A}$ denote the class of functions
$$\Phi : t\in\R \rightarrow \Phi(t) $$
where for almost every $t\in \R$, $\Phi(t)=\{\gamma_t^1,...,\gamma_t^N\}$ is a limit system of curves of class $\wpq{2}{p}$ and such that,  for almost every $\underline{t}, \overline{t}\in \R$, $\underline{t}< \overline{t}$, the following conditions are satisfied: 
\begin{itemize}
 \item[(i)]  $\Phi(\underline{t})$ and $\Phi(\overline{t})$ do not cross but may intersect tangentially;
 \item[(ii)] $\Int(\Phi(\overline{t})) \subseteq \Int(\Phi(\underline{t}))$ (pointwisely);
 \item[(iii)] if, for some $i$, $\mathcal{H}^1\left( (\gamma_{\overline{t}}^i) \setminus \overline{\Int(\Phi(\underline{t}))}\right) \neq 0$ then 
$$\mathcal{H}^1\left( [(\gamma_{\overline{t}}^i) \setminus \overline{\Int(\Phi(\underline{t}))}] \setminus (\Phi(\underline{t}))\right) = 0.$$ 
\end{itemize} 
\end{defi}

\boss\label{remarkiii}
 One may remark that, from condition $(ii)$ of Definition~\ref{A}, for every curve  $\gamma\in\Phi(\underline{t})$, $ (\gamma)  \cap  \Int(\Phi(\overline{t}))  =\emptyset.$ In fact if $x\in (\gamma)  \cap  \Int(\Phi(\overline{t}))$ then $x\in \Int(\Phi(\overline{t}))$ and $x\notin  \Int(\Phi(\underline{t}))$ which gives a contradiction with condition $(ii)$.
\eoss

\begin{defi}[The class $\mathscr{A}(u)$]\label{defiA}
Let $u\in \BV(\R^2)$. We define $\mathscr{A}(u)$ as the set of functions
$\Phi \in \mathscr{A} $ such that, for almost every $t\in \R$, we have 
$$(\Phi(t))\supseteq \partial^* \lt \quad \mbox{ (up to a $\mathcal{H}^1$-negligible  set)}$$
 and 
$$\lt = \Int(\Phi(t)) \quad \mbox{ (up to a ${\cal L}^2$- negligible  set)}.$$
In particular, if $u\in \cont^2(\R^2)$, we will denote as $\Phi[u]$ the function of  $\mathscr{A}(u)$ defined as
$$t \mapsto \partial \{u>t\}.$$
\end{defi}

In \cite{MN} we proved the following representation result for the relaxation problem for the Willmore functional on $\mathbb{R}^2$

\begin{thm}\label{principale}
Let $u\in \BV(\mathbb{R}^2)$ with  $\overline{F}(u,\mathbb{R}^2) < \infty$. Then $\ds\overline{F}(u,\mathbb{R}^2) = \underset{\Phi \in \mathscr{A}(u)} {\Min} G(\Phi).$
\end{thm}

The next proposition points out the relationship between the relaxation problem on $\mathbb{R}^2$ for a function with compact support and the relaxation problem on a suitable $\Omega\subset\R^2$:

\begin{prop}\label{legameY} Let $u\in \BV(\mathbb{R}^2)$ with compact support and such that $\overline{F}(u,\R^2) < \infty$. There exists an open bounded domain $\Omega$ such that
$$\overline{F}(u,\R^2)=\overline{F}_B^0(u,\Omega):=\inf\left\lbrace
\underset{h\rightarrow\infty}{\liminf}\;\ds\int_{\Omega}|\nabla u_h|(1+|\mdiv
\frac{\nabla u_h}{|\nabla u_h|}|^p)\,dx:\,\{u_h\}\in \cont^2_{\mathrm c}(\Omega), \;
u_h\overset{\lp{1}(\Omega)}{\longrightarrow}u\right\rbrace.$$
\end{prop}

As pointed out in \cite{MN} such a proposition is not true for the relaxation problem defined with $\cont^2$ instead of $\cont^2_{\mathrm c}$.

\section{Functionals defined on measures}\label{funzionalimisure}

Let $\mu,\nu$ be Radon measures on $\Omega\subset \rn$, $\mu$ positive, $\nu$
$\mathbb{R}^m$-valued and let $f:\mathbb{R}^m\rightarrow[0,\infty]$ be convex. We set 
$$G(\nu,\mu)=\int_{\Omega}f\left( \dfrac{\nu}{\mu}(x)\right)\;{d} \mu (x)
+\int_{\Omega}f^{\infty}\left( \dfrac{\nu^s}{|\nu^s|}(x)\right)\;{d} |\nu^s|(x)$$
where $\nu^s$ is the singular part of $\nu$ with respect to $\mu$ and
$f^{\infty}:\mathbb{R}^m\rightarrow\mathbb{R}\cup \{\infty\}$ is the recession
function of $f$ defined by
\begin{equation}\label{rec1}
 f^{\infty}(z)=\underset{t\rightarrow\infty}{\lim}\dfrac{f(z_0+tz)-f(z_0)}{t}
\end{equation}
where $z_0\in\mathbb{R}^m$ is any vector such that $f(z_0)<\infty$.

\par As stated in the theorem below, $G$ is lower semicontinuous under suitable assumptions.

\begin{thm}[\cite{AFP}, Thm 2.34]\label{sci}
Let $\Omega$ be an open subset of $\rn$ and $\nu,\nu_h$ be $\mathbb{R}^m$-valued
Radon measures on $\Omega$, $\mu,\mu_h$ positive Radon measures on $\Omega$. Let
$f:\mathbb{R}^m\rightarrow[0,\infty]$ be a convex lower semicontinuous function. If
$\nu_h\overset{*}{\rightharpoonup}\nu$ and $\mu_h\overset{*}{\rightharpoonup}\mu$ in
$\Omega$ then
$$G(\nu,\mu)\leq \underset{h\rightarrow \infty}{\liminf}\;G(\nu_h,\mu_h).$$
\end{thm}

Notice that if $f$ has superlinear  growth (i.e. $f^{\infty}(z)<\infty$ only if
$z=0$) then $G(\nu,\mu)<\infty$ only if $\nu<\!<\mu$ and in this case
$$G(\nu,\mu)=\int_{\Omega}f\left( \dfrac{\nu}{\mu}(x)\right)\;{d} \mu (x).$$
The next theorem will be useful in the sequel and  is a direct consequence of
Theorem \ref{sci}:

\begin{prop}[\cite{AFP}, Example 2.36]\label{236}
Let $\Omega$ be an open subset of $\rn$, $\nu, (\nu_h)$ $\mathbb{R}^m$-valued Radon
measures on $\Omega$, and $\mu,(\mu_h)$ positive Radon measures on $\Omega$. Let
$f:\mathbb{R}^m\rightarrow[0,\infty]$ be a convex lower semicontinuous function with
superlinear growth. If $\nu_h\overset{*}{\rightharpoonup}\nu$, $\mu_h\overset{*}{\rightharpoonup}\mu$ in
$\Omega$, $\nu_h<\!<\mu_h$ and
$\int_{\Omega}f(\nu_h/\mu_h)\;{d} \mu_h$
is bounded then $\nu<\!<\mu$ and 
$$G(\nu,\mu)\leq \underset{h\rightarrow \infty}{\liminf}\;G(\nu_h,\mu_h).$$
\end{prop}

%
%
\section{Varifolds}\label{teoriavarifold}
We collect below a few facts about varifolds. More details can be found in~\cite{Si,Allard}.
\subsection{Definitions}
We consider $G(n,k)$, $k\leq n$, the set of all $k$-dimensional subspaces of $\rn$
equipped with the metric 
$$\|S-T\|=\left(\sum_{i,j=1}^{n}| e_i\cdot P_S(e_j)-e_i\cdot P_T(e_j)|^2\right)^{1/2}\qquad\forall S,T\in G(n,k), $$
where $P_S$ and $P_T$ are the orthogonal projections of $\rn$ onto $S$ and
$T$, respectively, and $\{e_i\}_{i=1,\cdots,n}$ is the canonical basis of $\rn$.
$G(n,k)$ is called the  Grassmannian of all unoriented $k$-subspaces of $\rn$.
\par For a subset $\Omega$ of $\rn$ we define $G_{k}(\Omega)=\Omega\times G(n,k)$ equipped with the product metric. 

\begin{defi}[Varifolds] A $k$-varifold $V$ on $\Omega$ is a Radon measure on $G_{k}(\Omega)$. The 
weight measure  of V is the Radon measure on $\Omega$ defined by
$$\mu_{V}(U)=V(\pi^{-1}(U))$$
where $\pi$ is the projection $(x,S)\mapsto x$ of $G_{k}(\Omega)$ onto $\rn$.\\

\end{defi}

A very important class of varifolds is obtained from rectifiable sets.

\begin{defi}[Countably $\cont^r$-$k$-rectifiable sets]
$M\subseteq \mathbb{R}^{n}$ is a  countably $\cont^r$-$k$-rectifiable
set ($r\geq 1$) if
$$M = M_0 \cup \left(\bigcup_{i=1}^{+\infty} K_i\right)$$ where
$\mathcal{H}^{k}(M_0)=0$, $K_i\cap K_j=\emptyset$ if $i\neq j$ and for all $i\geq 1$
$K_i$ is a  subset of a $\cont^r$-$k$-manifold of $\mathbb{R}^{n}$.\\
$M$ is $\cont^r$-$k$-rectifiable if $M$ is countably $\cont^r$-$k$-rectifiable and
$\mathcal{H}^{k}(M)\leq+\infty$.\\
$M$ is (countably) $\cont^r$-$k$-rectifiable in $\Omega\subseteq \mathbb{R}^{n}$ if $M \cap
\Omega$ is (countably) $\cont^r$-$k$-rectifiable.
\end{defi}

Remark that if $M$ is a countably $\cont^r$-$k$-rectifiable set then for every  $x\in K_i$ 
 we can consider the tangent plan to  $K_i$ at $x$, denoted by $T_x$, and,  by the previous definition, the function 
$$x \mapsto T_x$$
is defined for $\mathcal{H}^{k}$-a.e. $x\in M$.

\begin{defi}[Rectifiable and integral varifolds]
$V$ is a rectifiable $k$-varifold if there exists a
$\cont^r$-$k$-rectifiable subset $M$ of $\Omega$ ($r\geq 1$)
such that
$$V=\theta\mathcal{H}^k\res M \otimes \delta_{T_x}$$ 
where $\theta $ is a
positive  $\mathcal{H}^k$-locally integrable function on $M$ called 
multiplicity of $V$. Then we denote  $V={\bf v}(M,\theta)$ and the weight measure of $V$ is
$$\mu_{V}(U)=\int_{M \cap U} \theta(x) \;{d} \mathcal{H}^k(x).$$
If $\theta(M)\subset\mathbb{N}$ then ${\bf v}(M,\theta)$ is called an integral varifold.\\
\end{defi}

\subsection{First variation and generalized mean curvature}

The first variation of a $k$-varifold $V$ is the functional given by
$$\delta V: \cont_{\mathrm c}(\Omega,\rn)\rightarrow \mathbb{R}$$
$$\delta V(X)=\int_{G_k(\Omega)}\mdiv_S X(x)\;{d} V(x,S),$$
where $\mdiv_S$ is the tangential divergence with respect to $S$.

\begin{defi}
A $k$-varifold is called a  Allard's varifold if it has locally bounded first
variation in $\Omega$, i.e. for each $W\subset \subset\Omega$ there exists a constant $0< C <\infty$ such that  
$$|\delta V(X)|\leq C\|X\|_{\lp{\infty}(W)}\mbox{\;\;\;}\forall X\in \cont_{\mathrm c}(W,\rn).$$ 
\end{defi}

If $V$ is a Allard's varifold then 
$\ds \|\delta V\|(W)=\sup \{|\delta V(X)|:  X\in \cont_{\mathrm c}(W,\rn), \|X\|_{\lp{\infty}(W)}\leq
1\}<\infty$ for each $W\subset\subset \Omega$ so, by the Riesz representation theorem,
$$\delta V(X)=-\int_{\Omega}\langle X, \nu\rangle \;{d} \|\delta V\|$$
where $\|\delta V\|$ is the total variation measure of $\delta V$ and $\nu$ is a
$\|\delta V\|$-measurable $\rn$-valued function with $|\nu|=1$ $\|\delta V\|$-a.e in
$\Omega$. By the Radon-Nikodym decomposition theorem,
$$\|\delta V\|=\dfrac{\|\delta V\|}{\mu_V}\mu_V+\sigma$$
where the derivative of $\|\delta V\|$ with respect to $\mu_V$ exists $\mu_V$-a.e.
and the measure $\sigma$, the singular part of $\|\delta V\|$ with respect to $\mu_V$ , is
supported on $Z$ such that
$$Z=\left\lbrace x\in \Omega : \dfrac{\|\delta V\|}{\mu_V}(x)=+\infty\right\rbrace \, ,\quad\quad \mu_V(Z)=0.$$

So, defining ${\bf H}_V(x)=\dfrac{\|\delta V\|}{\mu_V}(x)\nu(x)=-\ds\frac{\delta V}{\mu_V}(x)$, we can write
$$\delta V(X)=-\int_{\Omega}\langle X, {\bf H}_V\rangle \;{d} \mu_V - \int_{Z}
\langle X, \nu\rangle  \;{d} \sigma$$
for all $X\in \cont_{\mathrm c}(\Omega,\rn)$.

\begin{defi}
With the definitions above, ${\bf H}_V$ is the generalized mean curvature  of $V$, $Z$ the 
generalized boundary  of $V$, $\sigma$ the  generalized boundary measure of
$V$, and $\nu_{|Z}$ the  generalized unit conormal of $V$.
\end{defi}

\boss
Notice that if $V$ is a rectifiable $k$-varifold ${\bf
v}(M,1)$ associated with $M$ a $\cont^2$-manifold without boundary then, from  the divergence theorem for manifolds (see \cite{AFP}: Theorem 7.34), $\|\delta V\|<\!<\mu_V$ and the mean curvature  ${\bf H}_V$ coincides everywhere out of a $\mathcal{H}^k$-negligible set with the classical mean curvature of $M$.
\eoss

\boss{\bf (${\mathbf 2}$-Varifolds supported on ${\mathbf {\wpq{2}{p}}}$-curves)}\label{varW2p}
Let  $V={\bf v}(M,\theta)$ be the varifold on $\mathbb{R}^2$ associated with $M$ a closed curve in $\mathbb{R}^2$ of class $\wpq{2}{p}$ 
with $p> 1$, and with the density function $\theta$. $M$ admits a parametrization  (still denoted by  $M$) 
$M\in \wpq{2}{p}\left([0,L],\mathbb{R}^2 \right)$,
$$M(s)=(f(s), g(s))\;,\;\;f,g \in \wpq{2}{p}\left([0,L],\mathbb{R} \right)$$  
where $s$ is the arc-length parameter and $L$ is the length of the curve $M$. 
Then, by direct calculation, we will show that $\|\delta V\|<\!<\mu_V$ and the mean curvature of $V$  is a function of the weak second derivatives 
of $f$ and $g$. This fact can be generalized using Hutchinson's varifolds~\cite{HU1,HU2}.
\par Consider $X\in \cont_{\mathrm c}^1(\mathbb{R}^2, \mathbb{R}^2)$, $X(x) = (X^1(x), X^2(x))$, $\{e_1, e_2\}$ the canonical orthonormal basis of $\mathbb{R}^2$ and denote by $\langle \cdot, \cdot \rangle$ the usual scalar product in $\mathbb{R}^2$. Then 
$$
\begin{array}{ll}
\displaystyle{\mdiv_M X(M(s)) } & = \displaystyle{\langle e_1, \langle \nabla X^1 (M(s)), M'(s)\rangle M'(s)\rangle+\langle e_2, \langle \nabla X^2 (M(s)), M'(s)\rangle M'(s)\rangle =}\vspace{0.09cm}\\
 & \displaystyle{ = f'(s)\langle \nabla X^1 (M(s)), M'(s)\rangle + g'(s)\langle \nabla X^2 (M(s)), M'(s)\rangle}
\end{array} 
$$
and
$$
\begin{array}{ll}
\displaystyle{\delta V (X ) =  \int_M \theta \mdiv_M X \;{d} \mathcal{H}^1 } & \displaystyle{= \int_0^L \left[ f'(s)\langle \nabla X^1 (M(s)), M'(s)\rangle + g'(s)\langle \nabla X^2 (M(s)), M'(s)\rangle\right] \;{d} s}\\
& \displaystyle{= \int_0^L \left[ f'(s)\dfrac{d}{d s}\left[  X^1 (M(s))\right]  + g'(s)\dfrac{d}{d s} \left[  X^2 (M(s))\right] \right] \;{d} s}.\\
\end{array} 
$$
\par Now, integrating by parts and using the facts that  $X$ has compact support, $M$ is closed and  $f,g \in \wpq{2}{p}\left([0,L],\mathbb{R} \right)$, we get $\ds \delta V (X )= - \int_0^L \langle X(M(s)), (f''(s), g''(s))\rangle \;{d} s $, where $f'', g''$ are the weak second derivatives. It follows that $\ds\delta V(X) = -\int_M \langle X, {\bf H}_V \rangle \theta \;{d} \mathcal{H}^1$
where the curvature of varifold V is given by ${\bf H}_V (p)=M'' (M^{-1}(p))= (f''(M^{-1}(p)), g''(M^{-1}(p)))\;\; \forall \;p\in (M).$
Clearly $\|\delta V\|<\!<\mu_V$. By a  similar calculation we can generalize this remark to the varifolds $V={\bf v}(M, \theta)$ where $M$ is a system of curves of class $\wpq{2}{p}$ and $\theta$ the density function on $M$. 
\eoss

\section{Young measures}\label{Youngmeasure}

We collect below a few facts about Young measures, following \cite{KR}. 

\subsection{Definitions and general results}

Let $\Omega$ be a bounded Lipschitz domain of $\rn$ and let $f\in C(\Omega\times
\rn)$. By $\mathbb{B}^n$ we denote the open unit ball in $\rn$ and $\mathbb{S}^{n-1}=\partial \mathbb{B}^{n}$. We consider the following operator 
$$T: C(\Omega\times \rn) \rightarrow C(\Omega\times \mathbb{B}^n)$$ 
$$Tf(x,z):= (1-|z|)f\left( x, \dfrac{z}{1-|z|}\right)$$
and the property 

\begin{equation}\label{P}
\begin{array}{lll}
 \qquad  \mbox{T}f \;\mbox{extends into a bounded continuous function on}\;
\overline{\Omega\times \mathbb{B}^n}.
\end{array}
\end{equation}

We can define the Banach space $({\bf E}(\Omega;\rn), \|\cdot\|_{{\bf E}})$, where
$$\displaystyle{{\bf E}(\Omega;\rn)=\{f\in C(\Omega\times \rn) : f
\mbox{\;\;satisfies\;\;} \eqref{P}\}}$$
$$\|f\|_{{\bf E}}= \|\mbox{T}f\|_{\lp{\infty}(\overline{\Omega\times \mathbb{B}^n })}.$$
For example, a continuous function which is either uniformly bounded or  positively 1-homogeneous
 in its second argument (i.e. $f(x,sz)=sf(x,z)$, for all $s \geq 0$) belongs to ${\bf
E}(\Omega\times \rn)$.
Moreover, every $f\in{\bf E}(\Omega; \rn)$ has linear growth to infinity since 
$$|f(x,z)|= (1+|z|)Tf\left( x, \dfrac{z}{1+|z|}\right)\leq \|f\|_{\bf E}(1+|z|)
\qquad \mbox{   for all }x\in \Omega, z\in \rn.$$

For all $f\in{\bf E}(\Omega\times \rn)$ we define the  recession function
$f^{\infty} : \overline{\Omega}\times\mathbb{S}^{n-1}\rightarrow \mathbb{R}$ by
$$f^{\infty}(x,z) := \underset{\underset{\underset{t\rightarrow
\infty}{z'\rightarrow z}}{x'\rightarrow x}}{\lim}\dfrac{f(x', tz')}{t}.$$
Remark that for every convex function $f=f(z)$ belonging to ${\bf E}(\Omega; \rn)$  
 the previous definition coincides with \eqref{rec1} (this follows from continuity for convex functions and taking $z_0=0$ in \eqref{rec1}).

Before defining generalized Young measures it is convenient to recall some notations about parametrized measures.
For sets $E\subset \mathbb{R}^k$, $F\subset \mathbb{R}^l$ open or closed, a parametrized measure $(\nu_x)_{x\in E}$ is a mapping from $E$ to $\mathcal{M}(F)$, the set of Radon measures on $F$. It is said to be weakly* $\mu$-measurable, for some $\mu\in \mathcal{M}^+(E)$, if  the function $x\mapsto \nu_x(B)$ is $\mu$-measurable
for all Borel sets $B\subset F$. Here $\mu$-mesurability is the mesurability with respect to the $\mu$-completion of the
 Borel $\sigma$-algebra on $E$.
\par Let $\lp{\infty}_{w^*}(E, \mu, \mathcal{M}(F))$ denote the set of weakly* $\mu$-measurable parametrized measures $(\nu_x)_{x\in E}\subset \mathcal{M}(F)$ such that $\underset{x \in E}{\sup}\;|\nu_x|(F)<\infty$ (taking the essential supremum with respect to $\mu$). We will omit $\mu$ in the notation if it is the Lebesgue measure.

\begin{defi}[\cite{KR}] The set ${\bf Y}(\Omega,
\mathbb{R}^n)$  of all generalized Young measures is the set of all triplets
$(\nu_x,\lambda_{\nu},\nu_x^{\infty})$, simply written  $\nu$, such that : 
\begin{itemize}
\item[(i)] $\nu_x \in \lp{\infty}_{w^*}(\Omega, \mathcal{M}^1(\mathbb{R}^n))$ where the
map  $x\mapsto \nu_x $ is defined up to a $\mathcal{L}^n$-negligible set and with $x
\mapsto \langle\nu_x,|\cdot|\rangle \in \lp{1}(\Omega)$. $\nu_x$ is called 
oscillation measure.
\item[(ii)]$\lambda_{\nu}\in \mathcal{M}^+(\overline{\Omega})$. $\lambda_{\nu}$ is
called  concentration measure.
\item[(iii)] $\nu_x^{\infty}\in \lp{\infty}_{w^*}(\overline{\Omega},
\lambda_{\nu};\mathcal{M}^1(\mathbb{S}^{n-1}))$ where the map  $x\mapsto
\nu_x^{\infty} $ is defined up to a $\lambda_{\nu}$- negligible set. $\nu_x^{\infty}$
is called  concentration-angle measure.
\end{itemize}
\end{defi}

Therefore we can see ${\bf Y}(\Omega, \mathbb{R}^n)$ as a subset of ${\bf
E}(\Omega\times \rn)^{*}$ through the following duality pairing :
$$\langle\langle\nu,f\rangle\rangle:=
\int_{\Omega}\int_{\rn}f(x,z)\;{d} \nu_x(z)\;{{d} x} +
\int_{\overline{\Omega}}\int_{\mathbb{S}^{n-1}}f^{\infty}(x,z)\;{d} \nu_x^{\infty}(z)\;{d} \lambda_{\nu}(x)
$$

Then we can define the convergence for Young measures in the sense of duality:
\begin{defi}[{\bf Y}-convergence]
 A sequence $\{\nu_h\}\subset {\bf Y}(\Omega, \mathbb{R}^n)$ converges weakly* to
$\nu$ in ${\bf Y}(\Omega, \mathbb{R}^n)$, written $\nu_h \overset{\bf Y}{\rightarrow} \nu$,
if
$\langle\langle\nu_h,f\rangle\rangle\rightarrow\langle\langle\nu,f\rangle\rangle$
for all $f\in {\bf E}(\Omega\times \rn)$.
\end{defi}

Moreover, we have the following properties :

\begin{thm}[Closure, \cite{KR}, Cor. 1]\label{closureY}
 The set ${\bf Y}(\Omega, \mathbb{R}^n)$ is weakly* closed (as a subset of ${\bf
E}(\Omega\times \rn)^{*}$).
\end{thm}

\begin{thm}[Compactness, \cite{KR}, Cor. 2]\label{comp}
 Let $\{\nu_h\}\subset {\bf Y}(\Omega, \mathbb{R}^n)$ be a sequence  such that :\\
(i) the functions $x\mapsto \displaystyle{\int_{\rn}|\cdot|\;{d} \{\nu_h\}_x}$ are
uniformly bounded in $\lp{1}(\Omega)$;\\
(ii) the sequence $\{\lambda_{\nu_h}(\overline{\Omega})\}$ is uniformly bounded.\\
Then $\{\nu_h\}$ is weakly* sequentially relatively compact in ${\bf Y}(\Omega,
\mathbb{R}^n)$.
\end{thm}

Every Radon measure on $\overline{\Omega}$ can be associated with a Young measure:

\begin{defi}
Let $\mu \in \mathcal{M}(\overline{\Omega},\rn)$ with Radon-Nikodym
decomposition $\mu=\alpha\mathcal{L}^n\res\Omega+\mu_s$. The Young measure
$\nu_{\mu}$ associated with $\mu$ is defined by :
$$\nu_x=\delta_{\alpha(x)} \mbox{, \;\;  }\lambda_{\nu}=|\mu_s| \mbox{, \;\;  }
\nu_x^{\infty}=\delta_{\frac{\mu_s}{|\mu_s|}}$$

\end{defi}



Lastly, there exists a useful notion of  barycenter for Young measures:

\begin{defi}[{\bf Barycenter}]
The barycenter of $\nu \in {\bf Y}(\Omega, \mathbb{R}^n)$ is the measure ${\operatorname{Bar}}_\nu\in \mathcal{M}(\overline{\Omega},\rn)$ given by
$${\operatorname{Bar}}_{\nu}= \left(\int_{\rn} z \;{d} \nu_x \right)\mathcal{L}^n\res\Omega +
 \left(\int_{\mathbb{S}^{n-1}} z \;{d} \nu_x^{\infty} \right)\lambda_{\nu}.$$
\end{defi}

\subsection{Gradient Young measures}

\begin{defi}
 The  Young measure associated with $u\in \BV(\Omega)$ is the measure $\nu_{Du(x)}=(\nu_x,\lambda_{\nu},\nu_x^{\infty})$ with 
$$\nu_x=\delta_{\nabla u(x)} \mbox{, \;\;  }\lambda_{\nu}=|D^{s}u|\mbox{, \;\; 
} \nu_x^{\infty}=\delta_{\frac{D^{s}u}{|D^{s}u|}(x)}.$$
\end{defi}

Gradient Young measures can now be defined, see~\cite{KR}.

\begin{defi}
 We call $\nu\in {\bf Y}(\Omega, \mathbb{R}^n)$ a  gradient Young measure if
there exists a bounded sequence $\{u_h\}\subset \BV(\Omega)$ (called a {\it generating sequence}) such that
$\nu_{Du_h} \overset{\bf Y}{\rightarrow} \nu$. The set of gradient Young
measures is denoted as $\mbox{\bf {GY}}(\Omega, \mathbb{R}^n)$.
\end{defi}

Remark in particular that if $\nu$ is generated by  $\{u_h\}\subset \BV(\Omega)$, i.e. $\nu_{Du_h}
\overset{\bf Y}{\rightarrow} \nu$, then for all $f\in {\bf E}(\Omega\times \rn)$ we have 
$$
\begin{array}{rl}
\underset{h \rightarrow
\infty}{\lim}\langle\langle\nu_{Du_h},f\rangle\rangle 
 & =\displaystyle{\underset{h \rightarrow \infty}{\lim}\left[ \int_{\Omega}f(\nabla
u(x))\;{{d} x} + \int_{\overline{\Omega}}f^{\infty}\left(
x,\dfrac{D^{s}u}{|D^{s}u|}(x)\right) \;{d} |D^{s}u|(x)\right]= }\\
\displaystyle{=\langle\langle\nu,f\rangle\rangle} & \displaystyle{=
\int_{\Omega}\int_{\rn}f(x,z)\;{d} \nu_x(z)\;{{d} x} +
\int_{\overline{\Omega}}\int_{\mathbb{S}^{n-1}}f^{\infty}(x,z)\;{d} \nu_x^{\infty}(z)\;{d} \lambda_{\nu}(x)}.
\end{array}
$$

We will use a few connections shown in~\cite{KR} between general gradient Young
measures and $\BV$ functions:
\begin{prop}\label{stry}${ }$
\par \begin{enumerate}
\item Given $\nu\in \mbox{\bf {GY}}(\Omega, \mathbb{R}^n)$, all generating sequences $(u_h)\subset\BV(\Omega)$ converge weakly-* in $\BV$ to $u\in\BV(\Omega)$ such that $Du={\operatorname{Bar}}_\nu\res\Omega$ ($u$ is called an {\bf underlying deformation}).
\item  If $\{u_h\}\subset \BV(\Omega, \mathbb{R})$ is uniformly bounded in $\BV$, there exists a subsequence (not relabeled) such that $u_h\to u$ weakly-* in $\BV$ and $\nu_{Du_h} \overset{\bf Y}{\rightarrow} \nu$ for some $\nu\in \mbox{\bf
{GY}}(\Omega, \mathbb{R}^n)$ with $Du={\operatorname{Bar}}_\nu\res\Omega$. In general, $\nu$ may not coincide with $\nu_{Du}$.
\item If $u_h \rightarrow u$ strictly in $\BV(\Omega)$ then $\nu_{Du_h}
\overset{\bf Y}{\rightarrow} \nu_{Du}$.
\end{enumerate}
\end{prop}
The proof of $Du={\operatorname{Bar}}_\nu\res\Omega$ in 2. is easy but instructive. Testing the Young convergence with $f(x,z)=\langle g(x),z\rangle$ where $g\in
\cont^{\infty}_{\mathrm c}(\Omega, \rn)$ (thus $f^{\infty}=f$) yields
$$
\begin{array}{ll}
\underset{h \rightarrow \infty}{\lim}\langle\langle\nu_{Du_h},f\rangle\rangle & \displaystyle{= 
\underset{h \rightarrow \infty}{\lim}\left[ \int_{\Omega}\langle
g(x),\nabla u(x)\rangle\;{{d} x} + \int_{\Omega}\langle g(x),
\frac{D^{s}u}{|D^{s}u|}(x)\rangle \;{d} |D^{s}u|(x)\right]}\\
& \displaystyle{=\underset{h \rightarrow \infty}{\lim}\int_{\Omega}\langle
g(x),\;{d} Du_h(x)\rangle  }
\end{array}
$$
and then, taking the limit, $
\ds\langle\langle\nu,f\rangle\rangle = \int_{\Omega}\langle
g(x),\;{d} Du(x)\rangle.
$. Now, because of the choice of $f$, $\langle\langle\nu,f\rangle\rangle=\int_{\Omega}\langle
g(x),\;{d} {\operatorname{Bar}}_{\nu}\rangle
$ so 
$$
Du={\operatorname{Bar}}_{\nu}\res\Omega=\left(\int_{\rn} z \;{d} \nu_x
\right)\mathcal{L}^n\res\Omega +  \left(\int_{\mathbb{S}^{n-1}} z
\;{d} \nu_x^{\infty} \right)\lambda_{\nu}\res \Omega .
$$
The Radon-Nikodym decomposition of $\lambda_{\nu}$ with respect to $\mathcal{L}^n $ implies
\begin{equation}\label{decomp1}
\lambda_{\nu}\res \Omega=\frac{\lambda_{\nu}}{\mathcal{L}^n}(x)\mathcal{L}^n\res\Omega + \lambda_{\nu}^{s}\res\Omega
\end{equation}
therefore
\begin{equation}\label{nabla}
 \nabla{u}(x)=\int_{\rn} z \;{d} \nu_x 
+\frac{\lambda_{\nu}}{\mathcal{L}^n}(x)\int_{\mathbb{S}^{n-1}} z
\;{d} \nu_x^{\infty}\qquad\mathcal{L}^n\mbox{-a.e.}x\in \Omega
\end{equation}
$$
D^s u = \left( \int_{\mathbb{S}^{n-1}} z \;{d} \nu_x^{\infty}\right)
\lambda_{\nu}^{s}\res\Omega\qquad\mbox{and}\qquad \int_{\mathbb{S}^{n-1}} z \;{d} \nu_x^{\infty}\neq 0 \mbox{\;\;\;}|D^s
u|\mbox{-a.e.} x\in \Omega$$

\par We end this section with the Characterization Theorem for gradient Young measures, that we state in the specific form of the real-valued case, see~\cite[Thm 9]{KR} for the more general form, and~\cite[p.542]{KR} as well as~\cite[Remark 8]{KR} for a justification of the simplification.

\begin{thm}[Characterization, \cite{KR}, Thm 9] \label{caract}Let $\Omega\subset \rn$ be an open bounded Lipschitz domain. Then, a Young measure $\nu \in {\bf Y}(\Omega, \rn)$ satisfying 
$$\lambda_\nu(\partial \Omega) = 0$$
is a gradient Young measure, i.e. $\nu \in {\bf GY}(\Omega, \rn)$, if and only if
$$\ds \int_{\Omega}\int_{\rn}|z| \;{d} \nu_x(z)\;{{d} x}+\lambda_{\nu}(\Omega)< \infty $$
and there exists $u\in \BV(\Omega)$ such that $\operatorname{Bar}_\nu=Du$, i.e. $Du=\langle\operatorname{id},\nu_x\rangle{\cal L}^n\res\Om+\langle\operatorname{id},\nu_x^\infty\rangle\lambda_\nu$.
\end{thm}

\subsection{Identification of  gradient Young measures}\label{subsec:examples}

We recall, following~\cite{AB97,KR}, the classical techniques for the identification of a gradient Young measure $\nu$. If $\nu\in \mbox{\bf {GY}}(\Omega, \mathbb{R}^n)$, there exists a bounded
sequence $\{u_h\}\subset \BV(\Omega)$ such that $\nu_{Du_h} \overset{\bf
Y}{\rightarrow} \nu$.
\begin{itemize}
 \item[1)]\underline{identification of $\nu_x$}: test the Young convergence using 
$$f(x,z)= \Phi(x)\varphi(z)$$
 with $\Phi \in \lp{\infty}(\Omega)$ and $\varphi\in C_0(\rn)$. Then
$f^{\infty}=0$ (because $\varphi^{\infty}=0$) and we get 
 $$ \varphi(\nabla u_h) \rightharpoonup \int_{\rn}\varphi(z)\;{d} \nu_x(z)\mbox{\; in
\;}L^{1}(\Omega).$$
 Using this fact for all such $\varphi$, we can identify $\nu_x$. 
\par An important particular situation is when $\nabla u_h \rightarrow v$ $\mathcal{L}^n$-a.e. for some $v \in
L^{1}(\Omega)$ then, by the Dominated Convergence Theorem, we get $\varphi(\nabla u_h)
\rightarrow \varphi(v)$ in $L^{1}(\Omega)$ and therefore
 $$ \varphi(\nabla u_h) \rightharpoonup \varphi(v) \mbox{\; in \;}L^{1}(\Omega).$$
Therefore, if $\nabla u_h \rightarrow v$ a.e. then  $\nu_x=\delta_{v(x)}$ for $\mathcal{L}^n$-a.e. $x\in\Omega$.

 \item[2)]\underline{identification of $\lambda_{\nu}$ and $\nu_x^{\infty}$}: test the Young convergence using the function
$$f(x,z)=\Phi(x)|z|\varphi\left( x,\frac{z}{|z|}\right) $$
 with $\Phi \in C(\overline{\Omega})$ and 
$\varphi\in C(\mathbb{S}^{n-1})$. Then $f^{\infty}=\varphi$ and we have
$$\displaystyle{\int_{\Omega}\Phi(x)|\nabla u_h(x)|\varphi\left( \dfrac{\nabla u_h(x)}{|\nabla
u_h(x)|}\right) \;{{d} x} + \int_{\overline{\Omega}}\Phi(x)\varphi\left(
\dfrac{D^{s}u_h}{|D^{s}u_h|}(x)\right) \;{d} |D^{s}u_h|(x)} $$ 
$$\longrightarrow\int_{\Omega}\int_{\rn}\Phi(x)|z|\varphi\left( \dfrac{z}{|z|}\right)
\;{d} \nu_x(z)\;{{d} x} +
\int_{\overline{\Omega}}\int_{\mathbb{S}^{n-1}}\Phi(x)\varphi(z)\;{d} \nu_x^{\infty}(z)\;{d} \lambda_{\nu}(x)$$
The knowledge from 1) and testing with all such $\Phi$, $\varphi$ allows to identify $\lambda_{\nu}$ and $\nu_x^{\infty}$. 
\par In particular taking $\varphi=1$ so $f(x,z)=\Phi(x)|z|$, with $\Phi\in C(\overline{\Omega})$, we
get
$$\displaystyle{\int_{\Omega}|\nabla u_h(x)|\Phi(x) \;{{d} x} + \int_{\overline{\Omega}} \Phi(x) \;{d} |D^{s}u_h|(x) } 
\rightarrow\int_{\Omega}\int_{\rn}|z|\Phi(x)\;{{d} x} + \int_{\overline{\Omega}} \Phi(x) \lambda_{\nu}(x)$$
so 
\begin{equation}\label{lambda}
|Du_h| \overset{*}{\rightharpoonup} \left( \int_{\rn} |z|\;{d} \nu_x(z)\right)
\mathcal{L}^n\res\Omega +  \lambda_{\nu} \mbox{, \; in\;\;}
\mathcal{M}^+(\overline{\Omega}).
\end{equation}
\end{itemize}



We now illustrate on a few classical examples (see the one-dimensional counterparts in~\cite{AB97}) what kind of information can be carried by gradient Young measures. We will revisit later on these examples within the framework of Young varifolds.

\begin{ex}{\rm {\bf (Oscillations)}.\label{osc}
Let $n=2$, $\Omega=B(0,1)$ and 
 $$u_h(x)=\left\{
\begin{array}{ll}
 |x| - \frac{2k}{2^h} & \text{ if\;} |x|\in
\left[\frac{2k}{2^h},\frac{2k+1}{2^h}\right]  \vspace*{0.2cm}\\
 -|x| + \frac{2k+2}{2^h}&  \text{ if\;} |x|\in
\left[\frac{2k+1}{2^h},\frac{2k+2}{2^h}\right]  
\end{array} \right.\qquad \mbox{for }k=0,...,2^{h-1}-1$$

\begin{figure}[!h]
\begin{center}
\includegraphics[height=4cm]{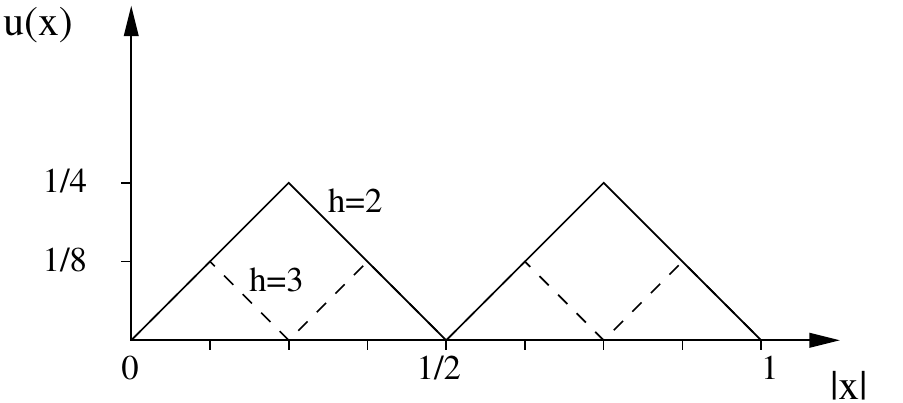}
\caption{A radial section of the graph of $u_h$}\label{cusp-ex-1}
\end{center}
\end{figure}

Since $\{u_h\}$ is uniformly bounded in $\BV(\Omega)$, extracting a subsequence (not relabeled) yields $\nu_{Du_h}\overset{\bf Y}{\rightarrow} \nu$
where $\nu_{Du_h}=(\nu_x^h,\lambda_\nu^h,\nu_x^{\infty,h})$ is defined by
$$\nu_x^h=\delta_{\nabla u_h(x)} \mbox{, \;\;  }\lambda_{\nu}^h=0 \mbox{,
\;\;  } \nu_x^{\infty,h}\text{ is arbitrary}.$$
Testing the Young convergence first with $f(x,z)=\Phi(x)\varphi(z)$, where $\Phi\in C(\overline{\Omega})$ and $\varphi\in C_0(\rn)$, and using polar coordinates and the Mean Value Theorem, then testing with $f(x,z)=\Phi(x)|z|$, it can be proved that
$$\nu_x=\dfrac{1}{2}\delta_{\frac{x}{|x|}}+\dfrac{1}{2}\delta_{-\frac{x}{|x|}}
\mbox{, \;\;  }\lambda_{\nu}=0,\quad\nu_x^\infty\text{ is arbitrary} .$$

}\end{ex}

\begin{ex}{\rm{\bf (Concentration)}.\label{conc}
Let $n=2$, $\Omega=B(0,2)$ and 
$$u_h(x)=\left\{
\begin{array}{ll}
 h(|x|-1)   &  \text{ if\;} |x|\in \left[1 ,1+\frac{1}{h}\right] \\
 h(1-|x|)   &  \text{ if\;} |x|\in \left[1-\frac{1}{h} ,1\right] \\
 0          &  \text{ otherwise}
\end{array} \right.$$

\begin{figure}[!h]
\begin{center}
\includegraphics[height=4cm]{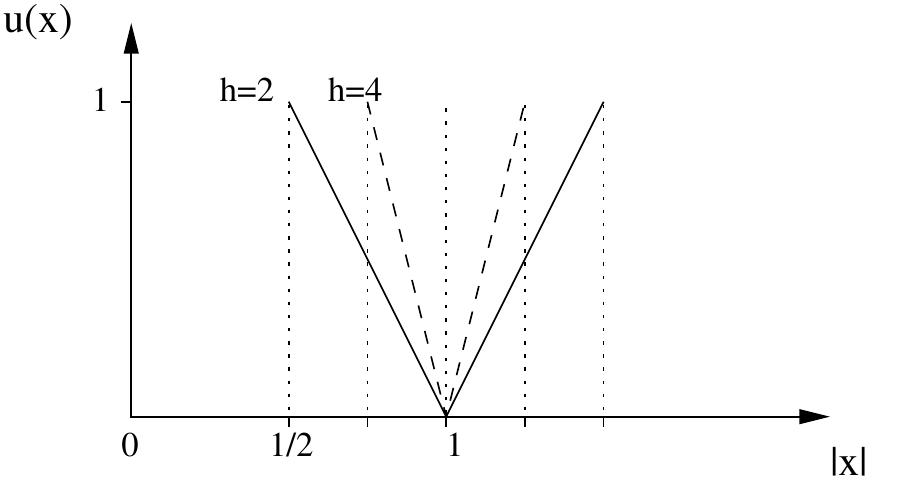}
\caption{A radial section of the graph of $u_h$}\label{cusp-ex-2}
\end{center}
\end{figure}

It is easily seen that $\{u_h\}$ is uniformly bounded in $\BV(\Omega)$ thus, possibly after extracting a
subsequence (not relabeled), $\nu_{Du_h}\overset{\bf Y}{\rightarrow} \nu$. Testing the Young convergence first with $f(x,z)=|z|$, then with $f(x,z)=\Phi(x)|z||\varphi(z/|z|)$  where $\Phi\in C(\overline{\Omega})$ and $\varphi\in C(\mathbb{S}^1)$ , it can be proved that
$$\nu_x=\delta_0 \mbox{, \;\;  }\lambda_{\nu}= 4\mathcal{H}^{1}\res
\partial{B}(0,1)\mbox{, \;\;  }
\nu_x^{\infty}=\dfrac{1}{2}\delta_{\frac{x}{|x|}}+\dfrac{1}{2}\delta_{-\frac{x}{|x|}}.$$

}\end{ex}

\begin{ex}{\rm {\bf (Diffuse concentration)}.\label{concdiff}
Let $n=2$, $\Omega=B(0,1)$ and 
$$u_h(x)=\left\{
\begin{array}{ll}
 h\left(|x|- \dfrac{k}{h}\right)  &    \text{ if\;} |x|\in \left[
\dfrac{k}{h},\dfrac{k}{h}+\dfrac{1}{2h^2}\right] \\
 h\left(\dfrac{k}{h}+\dfrac{1}{h^2}-|x|\right) &  \text{ if\;} |x|\in \left[
\dfrac{k}{h}+\dfrac{1}{2h^2},\dfrac{k}{h}+\dfrac{1}{h^2}\right]\\
 0     & \text{ otherwise} 
\end{array} \right.\qquad\mbox{where }k=0,...,h-1. $$

\begin{figure}[!h]
\begin{center}
\includegraphics[height=4cm]{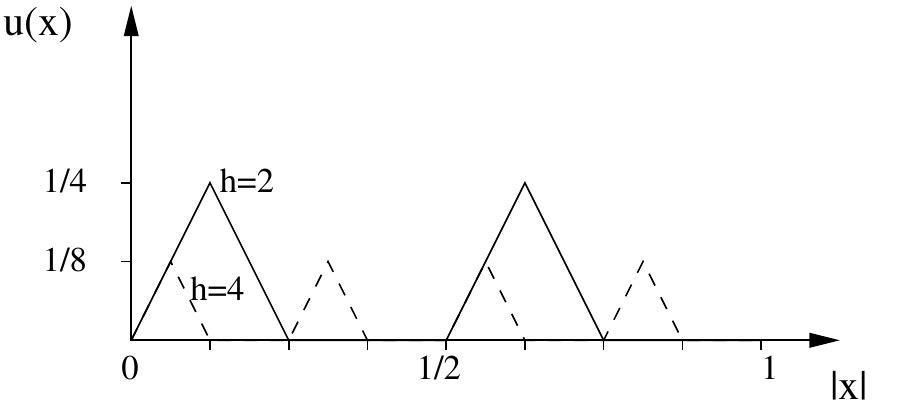}
\caption{A radial section of the graph of $u_h$}\label{cusp-ex-3}
\end{center}
\end{figure}

$\{u_h\}$ is clearly bounded in $\BV(\Omega)$ so, possibly after the extraction of a
subsequence (not relabeled), $\nu_{Du_h}\overset{\bf Y}{\rightarrow} \nu$. Testing the Young convergence first with $f(x,z)=|z|$, then with $f(x,z)=\Phi(x)|z||\varphi(z/|z|)$, where
$\Phi\in C(\overline{\Omega})$ and $\varphi\in C(\mathbb{S}^1)$, it can be proved that
$$\nu_x=\delta_0 \mbox{, \;\;  }\lambda_{\nu}= \mathcal{L}^2\res \Omega \mbox{,\;\;  }
\nu_x^{\infty}=\dfrac{1}{2}\delta_{\frac{x}{|x|}}+\dfrac{1}{2}\delta_{-\frac{x}{|x|}}.$$
This example shows that diffusion phenomena can also be generated by  sequences converging to zero.

\par More generally, all previous examples illustrate that the limit gradient Young measure is not determined  by the $\BV$ limit function but rather by the kind of sequence that generates it.}
\end{ex}

\section{Young varifolds}\label{Youngvarifold}

In this section, Young varifolds are defined and their basic properties are studied.

\begin{defi}[Young varifolds]
Let $\Omega$ be a bounded Lipschitz domain and let $V$ be a $(n-1)$-varifold. We say that $V$ is a  Young varifold if there exists  $\nu\in \mbox{\bf {Y}}(\Omega, \mathbb{R}^n)$ such that 
$$V(E\times A) = \int_{E}\int_{\{z\in \rn :  z^\perp\in A\}}|z|
\;{d} \nu_x(z)\;{{d} x}+ 
\int_{\overline{E}}\nu_x^{\infty}( \{ z\in \mathbb{S}^{n-1} : z^\perp \in A\} )\;{d} \lambda_{\nu} $$ 
for every $E\times A \subseteq G_{n-1}(\Omega)$, where $z^\perp$ is the element of $G(n,n-1)$ perpendicular to the space spanned by $z$.
Then $V$ is denoted by $V_\nu$ and is called  Young varifold associated with $\nu$. 
\par The weight measure is given by
$$\mu_{V_{\nu}}(E)=V_{\nu}(E\times G(n, n-1) )=\int_{E}\int_{\rn}|z|
\;{d} \nu_x(z)\;{{d} x}+\lambda_{\nu}(\overline{E}) \quad  \; \forall E\subseteq \Omega.$$
We denote by ${\bf {YV}}(\Omega, \mathbb{R}^n)$ the class of  Young varifolds.
\end{defi}

Remark that, for every $f\in \cont_{\mathrm c}(G_{n-1}(\Omega))$ the function $g(x,z)=|z|f\left( x, z^\perp  \right)$ belongs to ${\bf E}(\Omega, \rn)$. Moreover, as for every $s\in \mathbb{R}$ the linear spaces $(sz)^\perp$ and $z^\perp$ represent the same element of $G(n,n-1)$, $g$ is continuous and positively 1-homogeneous in $z$ and it is easy to check that $g^\infty = f$.  Thus we have $\ds\int_{G_{n-1}(\Omega)}f(x,S)\;{d} V_\nu(x,S)= \langle\langle \nu, |z|f( x, z^\perp  ) \rangle\rangle \qquad \forall f\in \cont_{\mathrm c}(G_{n-1}(\Omega)).$

\par As was mentioned in the introduction, $V_\nu$ may not be a rectifiable varifold since its projection on $\Omega$ might be a $n$-measure whereas the tangent measure lives in $G(n,n-1)$.
 
\par The following proposition shows that the  convergence of Young measures implies  the convergence of the associated Young varifolds:

\begin{prop}\label{yv}
If $\nu_{h} \overset{\bf Y}{\rightarrow} \nu$ then
$V_{\nu_{h}}\overset{*}{\rightharpoonup}V_{\nu}$.
\end{prop}

\dimo Given $f\in \cont_{\mathrm c}(G_{n-1}(\Omega))$, the Young convergence is tested with $g(x,z)=|z|f(x, z^\perp )\in {\bf E}(\Omega; \rn)$ .  Then if $\nu_{h} \overset{\bf
Y}{\rightarrow} \nu$ we have 
$$\int_{G_{n-1}(\Omega)}f \;{d} V_{\nu_{h}}=\langle\langle
g,\nu_{h}\rangle\rangle\rightarrow \langle\langle
g,\nu\rangle\rangle=\int_{G_{n-1}(\Omega)}f \;{d} V_{\nu}\;\;\;\;\;\forall f\in \cont_{\mathrm c}(G_{n-1}(\Omega))$$
so $\ds V_{\nu_{h}}\overset{*}{\rightharpoonup}V_{\nu}.$
\qed

Next proposition provides a sufficient condition for compactness in  ${\bf {YV}}(\Omega, \mathbb{R}^n)$. 

\begin{prop}[Compactness]\label{chiusura}
Let $\{V_h\}\subseteq {\bf {YV}}(\Omega, \mathbb{R}^n)$ be a sequence of Young varifolds such that $\underset{h}{\sup}\; \mu_{V_h}(\Omega) < \infty$. Then, possibly extracting a subsequence, $V_{h}\overset{*}{\rightharpoonup}V$ with $V\in {\bf {YV}}(\Omega, \mathbb{R}^n)$. 
\end{prop}

\dimo 
By definition of Young varifolds,  there exists a sequence of Young measures $\{\nu_h\}$ such that $V_h=V_{\nu_h}$ and by the uniform bound on $\mu_{V_h}$ we get 
$$\sup_h  \left[ \int_{\Omega}\int_{\rn}|z|
\;{d} (\nu_h)_x(z)\;{{d} x}+\lambda_{\nu_h}(\overline{\Omega}) \right]  <\infty.$$ 
Then, by Theorem \ref{comp} , there exists a  Young measure $\nu$ such that $\nu_h \overset{\bf Y}{\rightarrow}\nu$ (possibly extracting a subsequence) and  we get 
\begin{equation}\label{convi}
 \underset{h \rightarrow \infty}{\lim}\int_{G_{n-1}(\Omega)}f(x,S) \;{d} V_{\nu_{h}}(x,S)= \underset{h \rightarrow \infty}{\lim}\langle\langle
g,\nu_{h}\rangle\rangle= \langle\langle
g,\nu\rangle\rangle \qquad \forall f\in \cont_{\mathrm c}(G_{n-1}(\Omega))
\end{equation}
 where  $g(x,z)=|z|f(x, z^\perp )$. Then considering the Young varifold associated with $\nu$,
\eqref{convi} proves that $V_{h}\overset{*}{\rightharpoonup}V_{\nu}$ and the proposition ensues. 
\qed

\par The {\bf first variation of a Young varifold $V_\nu$}  is defined as
$$
\setlength{\extrarowheight}{0.4cm}
\begin{array}{llll}
 \delta V_\nu : &X\in\cont_{\mathrm c}(\Omega,\rn) &\longmapsto& \displaystyle{\int_{G_{n-1}(\Omega)} {\rm{div}}_S X(x)\; \;{d} V_{\nu}(x,S)}\\ 
&&&  \displaystyle{=\int_{\Omega}\int_{\rn}|z|{\rm{div}}_{ z^\perp}X\;\;{d} \nu_x(z)\;{{d} x}+\int_{\overline{\Omega}}\int_{\mathbb{S}^{n-1}}{\rm{div}}_{z^{\perp}}X\;\;{d} \nu_x^{\infty}(z)\;{d} \lambda_{\nu}.}
\end{array}
$$

\begin{ex}{\rm 
If $u\in \BV(\Omega)$ then  the Young varifold $V_{\nu_{Du}}$ associated with the gradient Young measure $ \nu_{Du} $ is defined as:
$$\int_{G_{n-1}(\Omega)}f(x,S) \;{d} V_{\nu_{Du}}(x,S)=\int_{\Omega}|\nabla u|f\left( x, \nabla u^\perp \right) \;{{d} x}+ 
\int_{\Omega}f\left(x, \frac{D^s u}{|D^s u|}^\perp  \right) \;{d} |D^s u|$$
for all $f\in \cont_{\mathrm c}(G_{n-1}(\Omega))$.
The weight measure is 
$$\mu_{V_{\nu}}=|Du|$$
and the first variation is
$$\delta V_{\nu_{Du}}(X)=\int_{\Omega}|\nabla
u|{\rm{div}}_{\nabla u^{\perp}}X\;{{d} x}+\int_{\Omega}{\rm{div}}_{{\frac{D^s
u}{|D^s u|}}^{\perp}}X\;{d} |D^s u|$$
for all $X\in \cont_{\mathrm c}(\Omega, \rn)$.
}\end{ex}

We can observe that, if $\{u_h\}$ is bounded in $\BV(\Omega)$, then, by Proposition
\ref{stry}, there exists a subsequence (not relabeled) such that $\nu_{Du_h}
\overset{\bf Y}{\rightarrow} \nu$ thus $V_{\nu_{Du_h}}\overset{*}{\rightharpoonup}V_{\nu}$. Furthermore, if
$u_h\rightarrow u$ strictly in $\BV$, Proposition~\ref{stry} implies that $V_{\nu_{Du_h}}\overset{*}{\rightharpoonup}V_{\nu_{Du}}.$

\boss({\bf Smooth functions})\label{reg}
If $u \in \cont^2(\Omega)$, $\Omega\subset\R^n$, then  for all $f\in \cont_{\mathrm c}(G_{n-1}(\Omega))$ we get
$$\int_{G_{n-1}(\Omega)}f(x,z) \;{d} V_{\nu_{Du}}(x,z)=\int_{\Omega}f(x,\nabla
u^\perp )|\nabla u|\;{{d} x}.$$
The weight measure is 
$\ds\mu_{V_{\nu_{Du}}}(A)=|Du|(A)=\int_{A}|\nabla u|\;{{d} x}$, $\forall A\subseteq \Omega .$ Moreover, for all $X\in \cont_{\mathrm c}(\Omega, \rn)$
$$\delta V_{\nu_{Du}}(X)=\int_{\Omega}|\nabla u|{\rm{div}}_{\nabla u^{\perp}}X\;{{d} x}.$$
\par Since $\nabla u$ is regular, the coarea formula yields $\ds\delta V_{\nu_{Du}}(X)=\int_{\mathbb{R}}\int_{\partial \{u>t\}\cap\Omega} {\rm{div}}_{\partial
\{u>t\}} X \;{d} \mathcal{H}^{n-1}\;{d} t$
because ${\nabla u}^{\perp}(x)$ is the tangent space at $x$ to the isolevel surface $\{y,\,u(y)=u(x)\}$. Moreover,
for a.e. $t\in \mathbb{R}$ the generalized mean curvature of the varifold $V_t={\bf v}(\partial \{u>t\}, 1)$ coincides with the mean curvature vector of the $\cont^2$ $(n-1)$-manifold $\partial \{u>t\}$. 
\par Now, a normal unit vector to $\partial \{u>t\}$ at $x\in \partial \{u>t\}$ is  $\nabla u(x)/|\nabla u(x)|$ and, denoting by $\mathbf{H}_{t}(x)$ the mean curvature (in the manifold sense !) of  $\partial \{u>t\}$ at $x$, we get
$$\mathbf{H}_{t}(x) = -({\rm{div}} \frac{\nabla u}{\left|\nabla u\right|}(x))\frac{\nabla u}{\left|\nabla u\right|}(x).$$
\par Then, using the representation formula for the first variation of rectifiable varifolds
and the coarea formula, it follows that
$$\delta V_{\nu_{Du}}(X)=-\int_{\mathbb{R}} \int_{\partial{\{u>t\}}\cap\Omega} \langle X,
\mathbf{H}_{t} \rangle \;{d} \mathcal{H}^{n-1}\;{d} t= \int_{\Omega}
|\nabla u| \langle X, ({\rm{div}} \frac{\nabla u}{\left|\nabla u\right|})\dfrac{\nabla u}{|\nabla u|}\rangle  \;{{d} x}.$$ 
\par Then the mean curvature vector of the varifold $V_{\nu_{Du}}$ is given by
$$\H_{V_{\nu_{Du}}}(x)= -\dfrac{\delta
V_{\nu_{Du}}}{\mu_{V_{\nu_{Du}}}}(x)=\mathbf{H}_t(x)=-({\rm{div}} \frac{\nabla u}{\left|\nabla u\right|}) \dfrac{\nabla u}{|\nabla u|}(x)
\mbox{\;\;\;}\mathcal{L}^n-a.e \mbox{\;in \;}\Omega.$$
Thus, for all $p> 1$, $\displaystyle{F(u, \Omega)=\int_{\Omega}\left|\nabla u\right|\left(1+\left|{\rm{div}}
\frac{\nabla u}{\left|\nabla
u\right|}\right|^p\right)\;{{d} x}= \int_{\Omega}\left[  1+ \left|\H_{V_{\nu_{Du}}}(x) \right|^p\right]  \;{d} \mu_{V_{\nu_{Du}}}}.$
\par This formula shows the relationship, in the case of a regular function, between  the generalized Willmore functional and the Young 
varifold associated with the Young measure $\nu_{Du}$ and it motivated our interest for the Willmore functional for Young varifolds studied in the next section.
\eoss

\section{The Willmore functional for Young varifolds}\label{Youngvarifold2}

In this section we extend the Willmore functional to  Young varifolds. We consider the class
$$\mbox{{\bf GY}}(u)=\{\nu \in\mbox{{\bf GY}}(\Omega, \rn) : {\operatorname{Bar}}_\nu\res \Omega=Du \}.$$

\begin{prop}\label{chiusura2}
 The set ${\bf GY}(u)$ is weakly* closed $($as a subset of $({\bf E}(\Omega, \rn))^*)$.
\end{prop}

\dimo Take $\nu$ from the weak closure of ${\bf GY}(u)$. By  Theorem \ref{closureY}, $\nu \in {\bf Y}(\Omega, \rn) $ and 
there exists a sequence $\{\nu_h\}\subset {\bf GY}(u)$ such that, for every $h \in \mathbb{N}$ and for every $f\in {\bf E}(\Omega, \rn)$,
$$\lvert\langle\langle\nu_h, f(x,z)\rangle\rangle-\langle\langle\nu, f(x,z)\rangle\rangle\rvert \leq \dfrac{1}{h} \quad\mbox{and}\quad \lvert\langle\langle\nu_h, |z| \rangle\rangle-\langle\langle\nu, |z|\rangle\rangle\rvert \leq \dfrac{1}{h}.$$
Now, for every $h$ $\nu_h\in{\bf GY}(u)$, and since $\wpq{1}{1}$-functions are dense in $\BV$ with respect to the strict convergence, it follows from Proposition \ref{stry} that for every $h\in\N$ there exists $u_h\in \wpq{1}{1}(\Omega)$ with 
$$\|u_h-u\|_{BV}\leq 1/h $$
and such that $\ds\lvert\int_\Omega f(x, \nabla u_h) d x - \langle\langle\nu_h, f(x,z)\rangle\rangle  \rvert \leq \dfrac{1}{h}$ and $\lvert\|\nabla u_h\|_{\lp{1}} - \langle\langle\nu_h, |z| \rangle\rangle\rvert\leq \dfrac{1}{h}.$ Being $\{u_h\}$ uniformly bounded in $\BV(\Omega)$, there exists a subsequence (not relabeled) such that
 $u_h \overset{*}{\rightharpoonup} u $ in $\BV$ and $\nu_{Du_h} \overset{{\bf Y}}{\rightarrow} \mu \in {\bf GY}(\Omega, \rn)$. 
Now, the two estimates above show  that necessarily  $\mu=\nu$  and we get ${\operatorname{Bar}}_\nu\res \Omega= {\operatorname{Bar}}_\mu\res \Omega= Du$,  hence $\nu \in {\bf GY}(u)$.
\qed

We can now define a class of Young varifolds that is suitable in the Willmore context.

\begin{defi}
The class $\mathbb{V}(u)$ of  Young varifolds associated with $\nu\in$ {\bf GY}$(u)$  is defined as
$$\mathbb{V}(u)=\{V_{\nu}\in {\bf {YV}}(\Omega, \mathbb{R}^n) : \nu\in{\bf
GY}(u)\, ,\; \|\delta V_\nu\|<\!<\mu_{V_\nu}\}.$$
\end{defi}

The following theorem shows a property of the weight measures of  Young varifolds that are associated with Young measures belonging  to ${\bf GY}(u)$.

\begin{prop}
For all $V\in \mathbb{V}(u)$, $\mu_{V_{\nu}}(\Omega)\geq |Du|(\Omega).$
\end{prop}

\dimo For every $\nu\in \mbox{{\bf GY}}(u)$, ${\operatorname{Bar}}_{\nu}\res \Omega= Du$ and, using the Radon-Nikodym decomposition
\eqref{decomp1}, 
$$\mu_{V_{\nu}}(\Omega)=\int_{\Omega}\int_{\rn}|z|
\;{d} \nu_x(z)\;{{d} x}+\lambda_{\nu}(\overline{\Omega})\geq \int_{\Omega}\left[ \int_{\rn}|z| \;{d} \nu_x(z)+\frac{\lambda_{\nu}}{\mathcal{L}^n}(x)\right] \;{{d} x} +
\lambda_{\nu}^s(\Omega)$$
Then, using  \eqref{nabla},  and reminding that $\nu_x^\infty$ is a probability measure, it follows that
$$
\begin{array}{ll}
\mu_{V_{\nu}}(A) & \geq 
\displaystyle{\int_{A}\lvert\int_{\rn} z \;{d} \nu_x +\frac{\lambda_{\nu}}{\mathcal{L}^n}(x)\int_{\mathbb{S}^{n-1}} z
\;{d} \nu_x^{\infty}\rvert d x +  \int_{A}\lvert \int_{\mathbb{S}^{n-1}} z \;{d} \nu_x^{\infty}\rvert
\lambda_{\nu}^{s}} \vspace*{0.2cm}\\
 &\displaystyle{ \geq \lvert\int_{A} \int_{\rn} z \;{d} \nu_x\;{d}x+\int_{A}\int_{\mathbb{S}^{n-1}} z \;{d} \nu_x^{\infty} \;{d}\lambda_{\nu}\rvert} =|Du(A)|
\end{array}
$$
for every  $A\subseteq \Omega$.
From the definition of total variation (see \cite{AFP}: Definition 1.4), we get 
$$\mu_{V_{\nu}}(\Omega)\geq |Du|(\Omega).$$
\qed

\par Remark that, since every  $u\in \BV(\Omega)$ can be approximated by a
sequence $\{u_h\}\subset \wpq{1}{1}(\Omega)$ strictly converging to $u$ in $\BV$, Proposition \ref{stry} implies that $\nu_{Du}\in {\bf GY}(u)$ 
for all $u\in \BV(\Omega)$. \par However  $V_{\nu_{Du}}\notin \mathbb{V}(u)$ in general (see Example \ref{non}) because 
it depends on the absolute continuity of $\|\delta\mu_{V_{Du}}\|$ with respect
to $\mu_{V_{Du}}$.

\begin{ex}{\rm ({\bf A case where $\mathbf{{V_{\nu_{Du}}\notin \mathbb{V}(u)}}$})\label{non}
We consider $E\subset \mathbb{R}^2$  like in Fig. \ref{nonasscont}, $\Omega$ an open set such that $E\subset\subset \Omega$ and
$u=\mathds{1}_{E}$. 
\begin{figure}[h]
\begin{center}
\includegraphics[height=3.5cm]{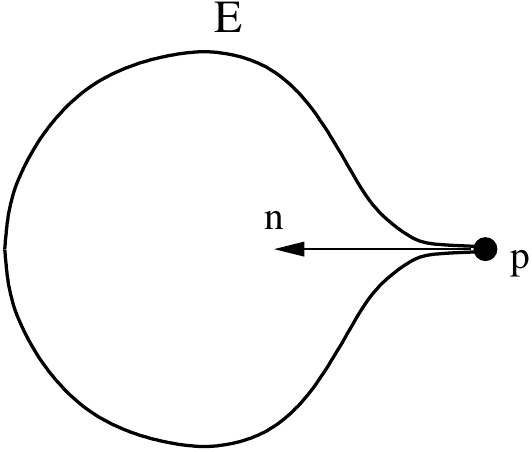}
\caption{$\delta V_{\nu_{Du}}$ has a singular component}\label{nonasscont}
\end{center}
\end{figure}

Remark that $V_{\nu_{Du}}=\mathbf{v}(\partial E, 1)$ and
$\mu_{V_{\nu_{Du}}}=\left|D^s u\right|$.
Nevertheless, using the  theory of rectifiable varifolds, it is easy to check that  $\left\|\delta{{V_{\nu_{Du}}}}\right\|$ is not absolutely continuous with respect to $\mu_{V_{\nu_{Du}}}$. In particular, denoting by $\sigma$ the generalized boundary measure of ${V_{\nu_{Du}}}$, we have $ \sigma = 2\bf{n}\delta_{p}$ where $\bf{n}$ is the unit vector drawn in Fig. \ref{nonasscont}.
}\end{ex}

The mean curvature vector and the Willmore functional associated with a Young varifold are defined as follows
\begin{defi}[Mean curvature of a Young varifold]
The generalized mean curvature vector of $V\in\mathbb{V}(u)$ is defined as the Radon-Nikodym derivative
$$ {\bf H}_V=-\frac{\delta V}{\mu_V}$$
\end{defi}
\begin{defi}[Willmore energy of a Young varifold]
The  Willmore energy of a Young varifold is defined as:
$$W:\mathbb{V}(u)\longrightarrow \mathbb{R} $$
$$W(V)=\int_{\Omega}\left( 1+ \left|{\bf H}_V\right|^p \right)
\;{d} \mu_{V} \quad p > 1.$$
\end{defi}

\par Remark that in general the class $\mathbb{V}(u)$ is not closed with respect to varifold
convergence because, given a sequence $(V_h)$ in $\mathbb{V}(u)$, the condition $\|\delta{V_h}\|<\!<\mu_{V_h}$ may not be preserved in the limit. A sufficient condition that ensures the preservation is a uniform bound (see Proposition~\ref{236}) 
$$\sup_h\int_{\Omega} \left| {\bf H}_{V_h}\right|^p  \;{{d} } \mu_{V_h}<+\infty.$$

\begin{oss}[Regular case]\label{wfr} 
{\rm If $u \in \cont^2(\Omega)$ and $F(u, \Omega)<\infty $ then $\nu_{Du}\in\mbox{{\bf
GY}}(u)$ and from Example \ref{reg} we have 
$\|\delta{V_{\nu_{Du}}}\|<\!<\mu_{V_{\nu_{Du}}}$. Then $V_{\nu_{Du}}\in \mathbb{V}(u)$ and $W(V_{\nu_{Du}})=F(u, \Omega)=\overline{F}(u, \Omega).$ Moreover, by the coarea formula,
$$\int_{G_{n-1}(\Omega)} f(x,S) \;{d} V_{\nu_{Du}}(x,S)= \int_{\mathbb{R}}\int_{G_{n-1}(\Omega)} f(x,S) \;{d} {{\bf v}(\partial\{u>t\}, 1)}{d} t$$
so, in the regular case,  the Young varifold $V_{\nu_{Du}}$ satisfies a slicing formula that involves the unit-density varifolds supported on the boundaries of the level sets of $u$.}
\end{oss}

We will now illustrate with the examples of Section~\ref{subsec:examples} a few situations where we can explicitly calculate the Willmore functional for Young varifolds and we will even show the continuity of the energy for the provided approximating sequences. These examples illustrate that Young varifolds are suitable for catching the limit energy.

\begin{ex}
{\rm Take the sequence of Young measures studied in Example~\ref{osc} and observe that, for all $h$, $\mu_{V_{\nu_h}}=\mu_{V_\nu}=\mathcal{L}^2 \res B(0,1)$. An easy calculation shows that for every $h$, $\H_{V_{\nu_h}}(x)=\H_{V_{\nu}}(x)=\dfrac{x}{|x|^2}$, $\forall x \in B(0,1)$,  and for $p\in(1,2)$,  
$$W(V_{\nu_h})=W(V_\nu)=\int_{B(0,1)}\left(1+\dfrac{1}{|x|^p}\right)\;{{d} x}.$$ 
}\end{ex}

\begin{ex}
{\rm Take now the sequence of Young measures studied in Example \ref{conc}. The limit varifold  satisfies $\mu_{V_\nu}=4\mathcal{H}^1 (\partial B(0,1))=8\pi$ and $\delta V_{\nu}(X)=-4\ds\int_{\partial B(0,1)}\langle X,x \rangle\;{d} \mathcal{H}^1$, therefore $\H_{V_{\nu}}(x)=x$, $\forall x \in \partial B(0,1)$, and $W(V_\nu)=16\pi.$
\par For  every $V_{\nu_h}$, $\mu_{V_{\nu_h}}=h\mathcal{L}^2\left(\Omega_h\right)  + \mathcal{H}^1(\partial\Omega_h)=8\pi$, where $\Omega_h=B(0,1+\dfrac{1}{h})\setminus B(0,1-\dfrac{1}{h})$ and it arises from an easy calculation that
$$W(V_{\nu_h})= 8\pi+ 2\pi h\int_{1-1/h}^{1+1/h} \dfrac{1}{r^{p-1}}\;{d} r + 2\pi
\left[\left(\dfrac{h}{h+1} \right)^{p-1} + \left(\dfrac{h}{h-1} \right)^{p-1}
\right] $$
By Lebesgue 's Theorem we  get  $W(V_{\nu_h})\rightarrow 8\pi+4\pi+4\pi=W(V_\nu).$
}\end{ex}

\begin{ex}{\rm \label{ex3}
Lastly, consider the sequence of Young measures studied in Example~\ref{concdiff}. For the Young varifold associated with  $\nu$ we have $\mu_{V_{\nu}}=\mathcal{L}^2\res \Omega$ and 
$
 \delta V_{\nu}(X)=\ds\int_{B(0,1)}{\rm{div}}_{{\frac{x}{|x|}}^{\perp}}X\;{{d} x}
$. An application of the coarea formula yields 
$\delta V_{\nu_h}(X) =-\ds\int_{B(0,1)}\langle X,\frac{x}{|x|^2}
\rangle\;{{d} x}$, therefore $\H_{V_{\nu}}(x)=\dfrac{x}{|x|^2}$ $\forall x \in B(0,1)$, and $\ds W(V_{\nu})=\pi+\int_{B(0,1)}\frac{\;{{d} x}}{|x|^{p}}==\dfrac{4-p}{2-p}\pi.$ Moreover, 
$$\ds W(V_{\nu_{h}})=\pi\sum_{k=0}^{h-1}\dfrac{2kh+1}{h^3}+2\pi\sum_{k=0}^{h-1} h \int_{\frac{k}{h}}^{\frac{k}{h}+\frac{1}{h^2}}\dfrac{\;{d} r}{r^{p-1}}\rightarrow\pi+2\pi\int_{0}^{1}\dfrac{\;{d} r}{r^{p-1}}=\dfrac{4-p}{2-p}\pi,$$
therefore  $W(V_{\nu_{h}})\rightarrow \dfrac{4-p}{2-p}\pi=W(V_{\nu}).$
}\end{ex}

\boss\label{rex3}
 Remark that, in all these examples, we have convergence of the Willmore energy, but the limit varifold  is not the varifold associated with $Du$, 
i.e. $V_{\nu}\neq V_{\nu_{Du}}$. In fact the sequence $\{u_h\}$ converges to $u=0$ weakly* in $\BV$ but $\nu_{Du_h}$ 
does not converge (in the sense of Young measures) to the gradient Young measure associated with $0$.
For instance, in Example~\ref{ex3}, the sequence of
gradients creates some curvature at the limit that is captured by the diffuse part
$\lambda_{\nu}$.
\par Moreover, since $u$ is identically $0$, its level sets are empty for all positive levels but $\mu_{V_{\nu}}\neq 0$. This means that  the level lines
of $u$ {\it do not provide any} information about the Young measure  generated by $\{u_h\}$ so it will not be possible in general to  write a slicing formula (like in Remark \ref{wfr}) that links the measures belonging to ${\bf GY}(u)$ and the level sets of $u$.
\eoss

\section{Relaxation of the generalized Willmore functional and Young varifolds}\label{rappresentazione}

This section is devoted to the minimum problem associated with the Willmore functional for  Young varifolds and its relationship with the relaxation problem for $F$. In the remaining, $\Omega$ is an open, bounded Lipschitz domain.
\begin{thm}\label{minVu}
Let $u\in \BV(\Omega)$ with $\overline{F}(u,\Omega) < \infty$. 
Then
\begin{enumerate}
\item $\mathbb{V}(u) \neq \emptyset$
\item The problem $ \Min \;\{W(V):V\in\mathbb{V}(u)\}$ has a solution, 
\item $\overline{F}(u,\Omega) \geq \Min \{W(V):V\in\mathbb{V}(u)\}.$
\end{enumerate}
\end{thm}

\dimo Let $\{u_h\}\subset \cont^2(\Omega)$ be a sequence converging to $u$ in $\lp{1}(\Omega)$ and such that $\overline{F}(u, \Omega)= \underset{h \rightarrow \infty}{\lim} F(u_h, \Omega)$. We can also suppose that $F(u_h, \Omega)$ is uniformly bounded so, possibly taking a subsequence, $u_h \overset{*}{\rightharpoonup}u$ in $\BV(\Omega)$.
Then, by Remark \ref{reg}, the Young varifolds $V_{v_{Du_h}}$ associated with the gradient Young measures ${v_{Du_h}}$ satisfy
$$F(u_h, \Omega) =W(V_{v_{Du_h}}), \quad \forall h.$$
Moreover, by Proposition \ref{stry} and possibly taking a subsequence, we can assume that $\nu_{Du_h} \overset{{\bf Y}}{\rightarrow} \tilde\nu  \in {\bf GY}(\Omega, \mathbb{R}^n)$ and hence $V_{\nu_{Du_h}} \overset{*}{\rightharpoonup} V_{\tilde\nu}$. Clearly, as  $u_h \overset{*}{\rightharpoonup}u$ in $\BV(\Omega)$, we have ${\operatorname{Bar}}_{\tilde\nu} \res\Omega = Du$ which implies that $\tilde\nu \in {\bf GY}(u)$.
\par Then, we deduce from Proposition~\ref{236} that $\| \delta V_{\tilde\nu} \| <\!< \mu_{V_{\tilde\nu}}$ and
$$ W(V_{\tilde\nu}) \leq \underset{h \rightarrow \infty}{\liminf}\; W(V_{\nu_{Du_h}}),$$
Therefore $V_{\tilde\nu} \in \mathbb{V}(u)$  and,
\begin{equation}
\overline{F}(u,\Omega) \geq \inf \{W(V):V\in\mathbb{V}(u)\}\label{fplusgrand}
\end{equation}

Let $\{V_h\} \subset \mathbb{V}(u)$ be a minimizing sequence such that $W(V_h)$
is uniformly bounded, therefore $\{\mu_{V_h}(\Omega)\}$ is uniformly bounded. By Theorem~\ref{comp}  and Proposition~\ref{chiusura2}, there exists  $\nu \in {\bf GY}(u)$ such that (possibly extracting a subsequence) $\nu_h \overset{\bf{Y}}{\rightarrow}\nu$ and, by Proposition \ref{yv}, $V_h \overset{*}{\rightharpoonup} V_\nu $.
\par
\begin{figure}[!htp]
\begin{center}
\includegraphics[width=4cm]{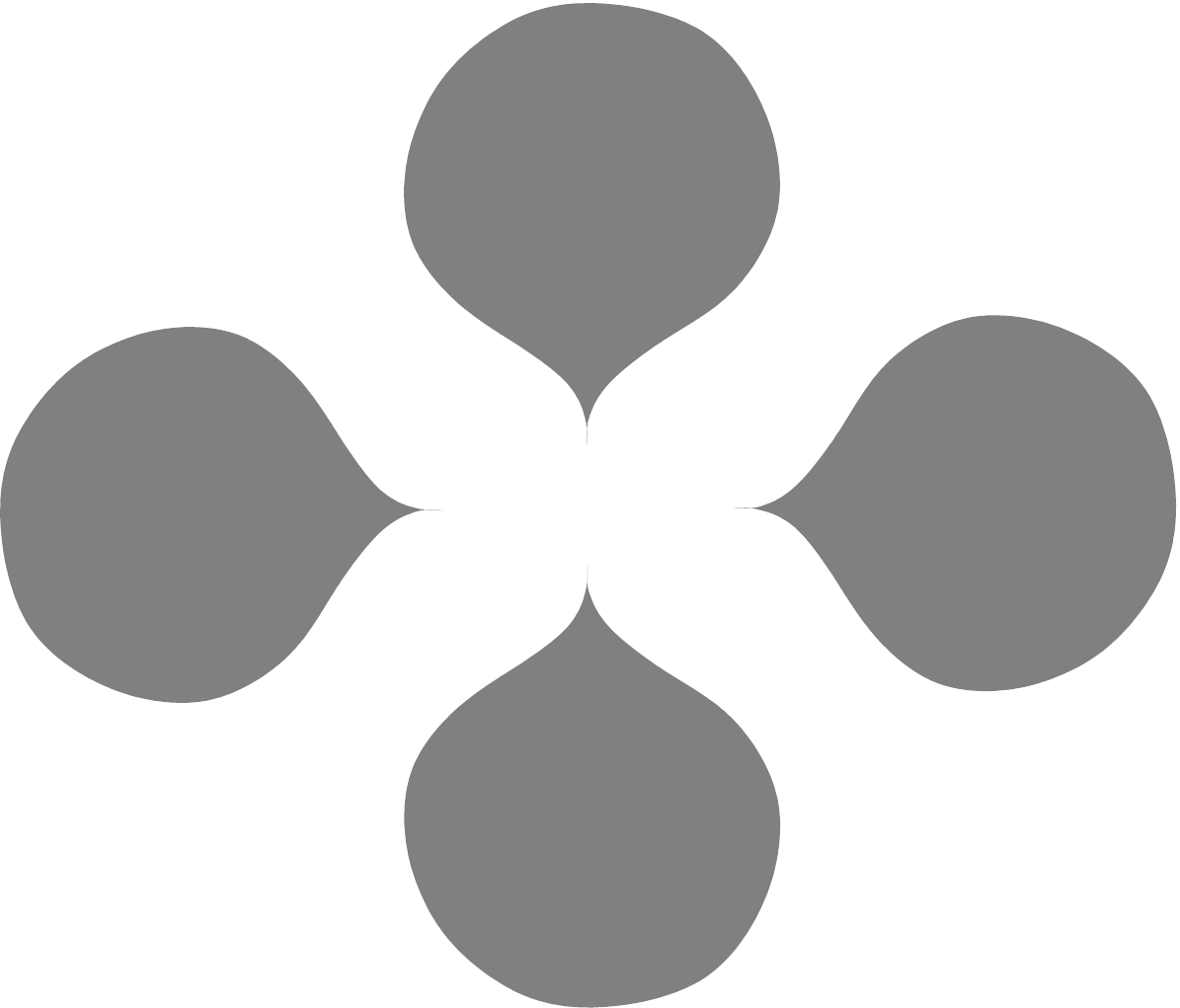}\qquad\includegraphics[width=4cm]{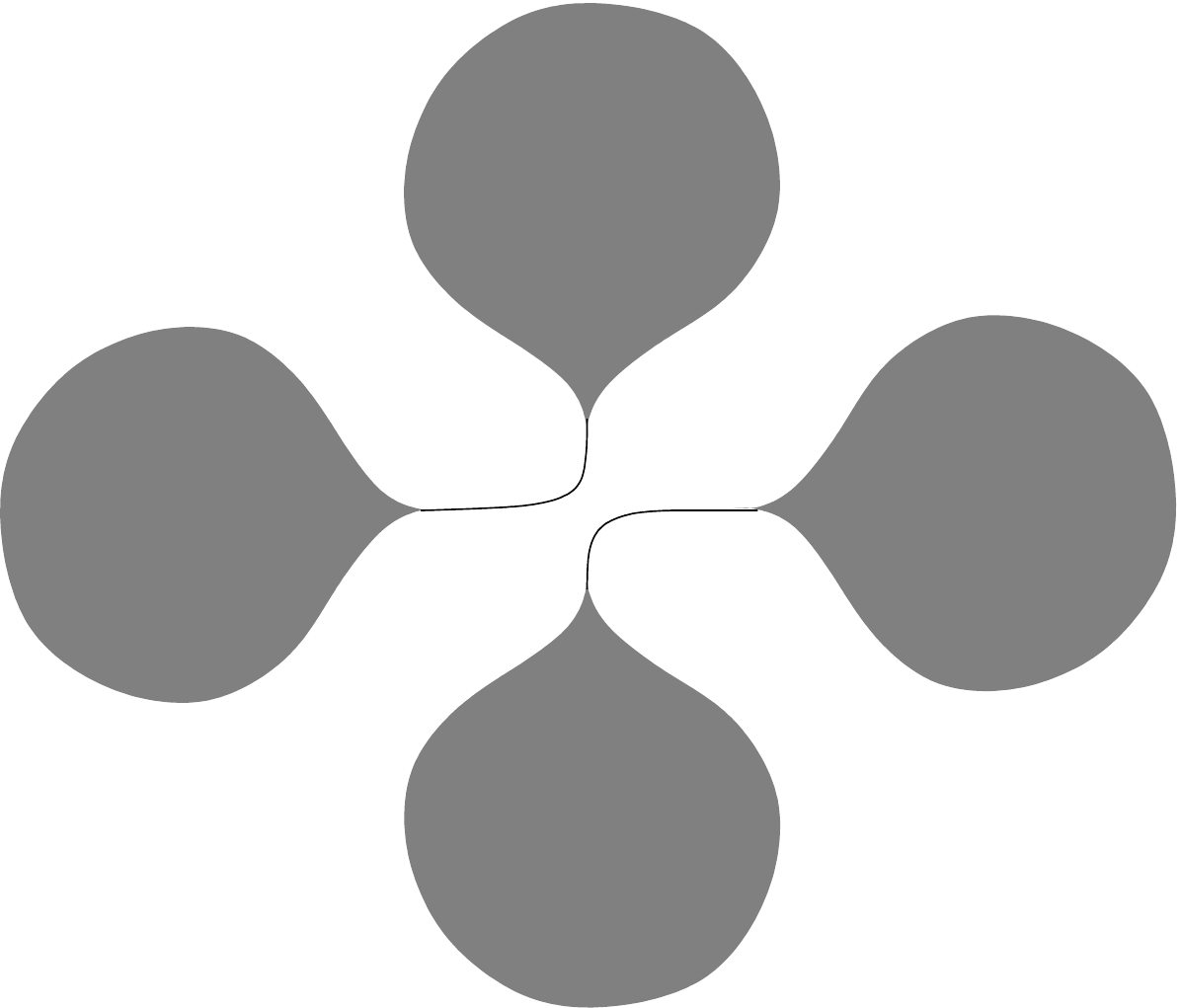}\qquad\includegraphics[width=4cm]{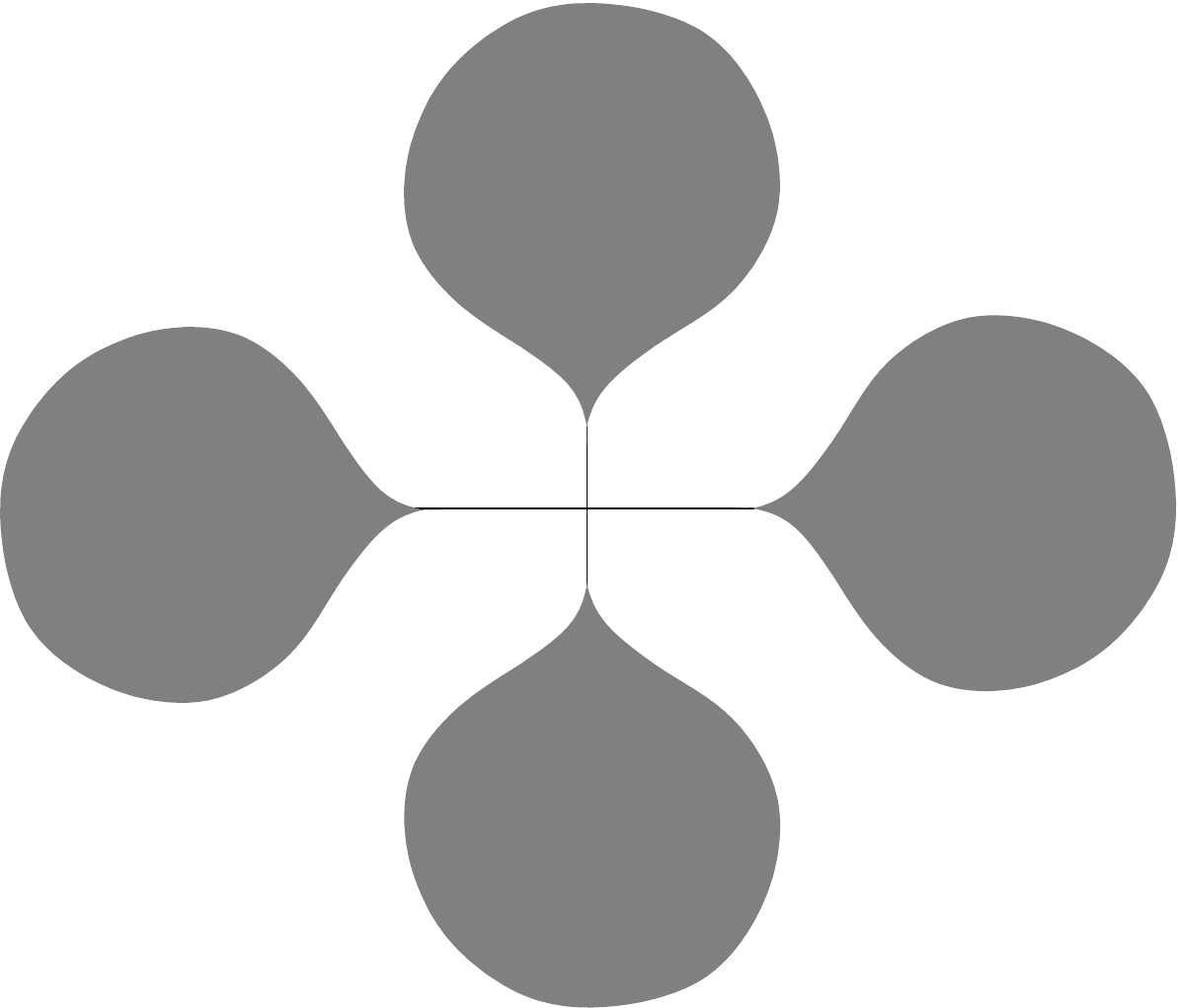}
\caption{Left : a set $E$. Middle: a limit configuration showing that $\oF(\one{E})<\infty$. Right: this configuration shows the support of a Young varifold $V\in{\mathbb V}(\one{E})$ such that $W(V)< \oF(\one{E})$.}\label{fig:cex}
\end{center}
\end{figure}

 Moreover, by Proposition~\ref{236},  $\| \delta V_{\nu} \| <\!< \mu_{V_\nu}$, $V_\nu \in \mathbb{V}(u)$  and 
$$W(V_\nu)\leq \underset{n\rightarrow\infty}{\liminf}\; W(V_h)$$
which proves that 
$$W(V_\nu)=\mbox{Min}\;\{W(V):V\in\mathbb{V}(u)\}.$$
Lastly, 3. follows from~\eqref{fplusgrand}.
\qed
As a by-product of the previous proof we have the 
\begin{coro}\label{corocomp}
For all $A>0$, the set $\mathbb{V}_A(u)=\{V\in\mathbb{V}(u),\; W(V)\leq A\}$ is sequentially compact, and $W$ is lower semicontinuous on it.
\end{coro}
\boss
Three questions arise naturally from Theorem~\ref{minVu}:
\begin{enumerate}
\item[Q.1] Does the equality hold in Theorem~\ref{minVu}, 3. ?
\item[Q.2] If not, what additional assumption should be taken to guarantee it ?
\item[Q.3] Can we at least always find $V\in\mathbb{V}(u)$ such that $\overline{F}(u,\Omega) =W(V_\nu)$?
\end{enumerate}
The answer to the first question is {\it negative}:
\begin{prop}\label{propcex}
There exists $u\in \BV(\Omega)$ such that $\overline{F}(u,\Omega) < \infty$ and
$$\overline{F}(u,\Omega) > \Min \{W(V):V\in\mathbb{V}(u)\}.$$
\end{prop}
\dimo
The counterexample is the function $u=\one{E}$ shown in Figure~\ref{fig:cex}, left.
The middle figure shows the limit configuration obtained with a smooth sequence $(u_h)$ converging weakly to $\one{E}$ in $\BV$ and such that $\sup_h F(u_h)<+\infty$. If the assumption $\sup_h F(u_h)<+\infty$ is dropped then the right configuration can be obtained with a suitable smooth sequence $(w_h)$ (just taking the previous situation and allowing the creation of a double right angle at the center). Associating $w_h$ with a Young measure $\nu_h$, and passing to the limit as Young measures, yields a limit Young measure $\nu$ whose support is the topological boundary of the right figure. An easy calculus shows that the associated Young varifold $V_\nu$ has no curvature at the center, therefore the central cross has no energy, thus $W(V_\nu)<\oF(\one{E})$.

\par Interestingly, the same counterexample can be used to show that the functionals 
$$ \frac{1}{2c_0}\int_{\R^n}\left( \epsilon |\nabla u_\epsilon|^2 + \frac{1}{\epsilon}  \Psi(u_\epsilon) \right) dx+\frac{1}{4c_0\epsilon} \int_{\R^n} \left( 2\epsilon \Delta u_{\epsilon} - \frac{1}{\epsilon} \Psi'(u_{\epsilon}) \right)^2 dx,$$
(where $\Psi(s)=(1-s^2)^2$ and $c_0=\int_{-1}^1\sqrt\Psi(t)dt$) 
do not $\Gamma$-converge to $\oF(\one{E})$ as $\epsilon\to 0$ (but the $\Gamma$-convergence holds in dimensions $2$ and $3$ when $E$ is smooth). The reason why such singular structure can be obtained is due to the existence of solutions of the Allen-Cahn equation with singular nodal set.

\par Clearly, the class of Young varifolds is richer than what is strictly needed to ensure the equality in Theorem~\ref{minVu}, 3. An even more extreme example is provided by the set $E$ in Figure~\ref{fig:triple} below: it follows from the results in~\cite{BDP} that $\oF(\one{E})=+\infty$ (because there is are oddly many cusps), yet $\one{E}$ can be approximated by a sequence of smooth functions with uniformly bounded $\BV$ norm, and the associated Young varifolds converge to a limit Young varifold $V_\nu$ whose energy near the triple point is null because it has no singular part and null curvature (being the three angles equal, the singularities compensate). It follows that $W(V_\nu)<+\infty$. Again, the framework of Young varifolds is suitable for catching the limit configuration.
\begin{figure}[!htp]
\begin{center}
\includegraphics[width=4cm]{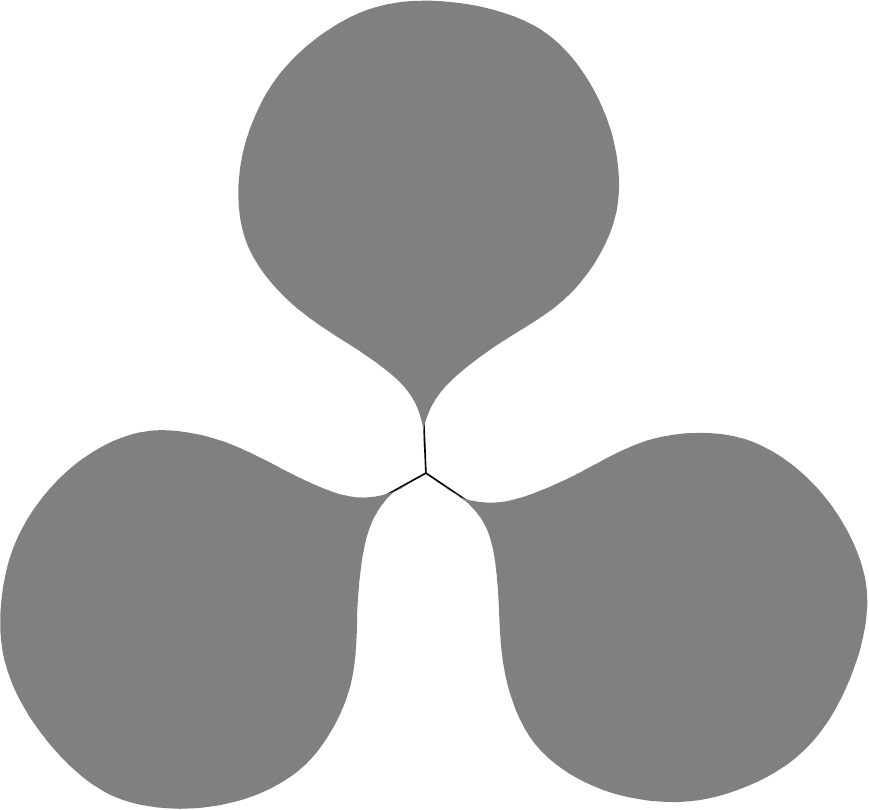}
\caption{A set $E$ such that $\oF(\one{E})=\infty$ but that can be associated with a Young varifold having finite energy.}\label{fig:triple}
\end{center}
\end{figure}
\qed

Being the answer to question Q.1 negative, what about Q.2?  To avoid the singular situations of Figures~\ref{fig:cex}, ~\ref{fig:triple}, it is enough to restrict the class of Young varifolds to those that can be approximated by sequences of smooth functions uniformly controlled in energy, i.e. working with the class
$$\mathbb{V}_B(u)=\{V_{\nu}\in {\bf {YV}}(\Omega, \mathbb{R}^n) : \nu\in{\bf
GY}(u)\, ,\; \|\delta V_\nu\|<\!<\mu_{V_\nu},\;\exists (u_h)\overset{*}{\rightharpoonup}u\mbox{ in }\BV,\,F(u_h)\leq B\}.$$
We do not know however if using such restriction is enough to guarantee the equality in Theorem~\ref{minVu}, 3.  To study in~\cite{BM} the representation by varifolds in dimension $2$ of the relaxed elastica functional, Bellettini and Mugnai need varifolds having a unique tangent line at {\it every} point. In our context and for our more general functional, one can reasonably think that a more intrinsic condition for the equality to hold could also involve tangential conditions, i.e. the existence of tangent planes along the support of $\lambda_\nu$, together with orthogonality conditions on $\nu_x^\infty$.  This is however a completely open problem.

\medskip
Let us end this (long!) remark with a few words on question Q.3. The answer is {\bf positive} in dimension $2$, as stated in Theorem~\ref{sistemaYoung} below where we use the connection beween systems of curves belonging to $\mathscr{A}(u)$ and  Young varifolds. In particular, given a system of curves, one can define a Young varifold with a coarea structure similar to the one in Remark~\ref{wfr}. The proof uses the following theorem on the locality of curvature for integral $1$-varifolds.
\eoss
\begin{thm}[Locality of mean curvature, \cite{LM}, Thm 2.1]\label{localH}
Let $V_1={\bf v}(M_1, \theta_1)$, $V_2={\bf v}(M_2, \theta_2)$ be two integral $1$-varifolds in $\mathbb{R}^{n}$ 
and let $\mathbf{H}_1, \mathbf{H}_2$ be their generalized curvature vectors.
Then
$$\mathbf{H}_1=\mathbf{H}_2\mbox{\;\;}\mathcal{H}^{1}\mbox{-a.e. on } M_1\cap M_2.$$ 
\end{thm}

\begin{thm}[Representation of the relaxation in dimension ${\mathbf 2}$]\label{sistemaYoung}
Let $u\in \SBV(\mathbb{R}^2)$ with compact support and $\oF(u)<\infty$. Let $\Phi\in \mathscr{A}(u)$ with $G(\Phi)<\infty$ and suppose that there exists a bounded Lipschitz domain $\Omega$ such that 
\begin{equation}\label{ipotbordo}
  \bigcup_{t\in \mathbb{R}}(\Phi(t))\subset\subset \Omega. 
\end{equation}
Then, there exists a Young measure  $\nu\in {\bf GY}(\Omega, \mathbb{R}^2)$ such that  $V_\nu \in \mathbb{V}(u)$ and
$$G(\Phi) = W(V_\nu).$$
In particular, there exists $\nu\in {\bf GY}(\Omega, \mathbb{R}^2)$ such that $V_\nu \in \mathbb{V}(u)$ and
$$\oF(u)=W(V_\nu)$$
\end{thm}

\dimo  We first remark that, since   $\Phi(t)$ is a system of curves of class $\wpq{2}{p}$ for a.e. $t$, 
the 1-rectifiable varifold $V_t={\bf v}(\Phi(t),\theta_{\Phi(t)})$ ($\theta_{\Phi(t)}$ is the density of the system $\Phi(t)$, 
see Definition \ref{sistemi}) is such that $$\|\delta V_t\|<\!<\mu_{V_t}$$ (see Remark \ref{varW2p}) so for every $X \in \cont_{\mathrm c}(\Omega, \mathbb{R}^2)$ we have
\begin{equation}\label{azionelivelli}
  \delta V_t(X) = - \int_{\Omega}\langle X, {\bf H}_{V_t}\rangle \;{d}  \mu_{V_t}
\end{equation}
where  ${\bf H}_{V_t}$ is the mean curvature of the varifold ${\bf v}(\Phi(t),\theta_{\Phi(t)} )$.
\par  Let $\Phi\in\mathscr{A}(u)$  and consider the varifold $V$, supported on $\bigcup_{t\in \mathbb{R}} (\Phi(t))\times G(2,1)$, defined as
\begin{equation}\label{varifoldfi}
 V(E \times A) = \int_{\mathbb{R}} {\bf v}(\Phi(t),\theta_{\Phi(t)})(E \times A) \;{d} t \;,\;\; \forall \;E \times A \subset G_1(\Omega)
\end{equation}
where $\theta_{\Phi(t)}$  is the density of the system $\Phi(t)$.
\par Then for every $f\in \cont_{\mathrm c}(G_{1}(\Omega))$ we get 
\begin{equation}\label{azionemiovar}
 \int_{G_{1}(\Omega)} f(x, S) \;{d}  V(x,S) = \int_{\mathbb{R}} \int_{\Phi(t)} \theta_{\Phi(t)}(x) f(x, n_t(x)^\perp) \;{d} \mathcal{H}^1 \;{d} t 
\end{equation}
where $n_t(x)$ denotes a unit normal to the system $\Phi(t)$ at $x$ and $\theta_t(x)$ denotes the density of the system $\Phi(t)$ at $x$. 
\par By  definition of $\mathscr{A}(u)$, for $\mathcal{L}^1$-a.e. $t$ $(\Phi(t)) \supseteq \partial^* \{u> t\}$ (up to a $\mathcal{H}^1$-negligible set) so
\begin{equation}\label{primadecomposizione}
\begin{array}{l}
\displaystyle{\int_{G_{1}(\Omega)} f(x, S) \;{d}  V(x,S) =\int_{\mathbb{R}} \int_{\partial^*\{u>t\}} f(x, n_t(x)^\perp)\;{d} \mathcal{H}^1 {d}t\,\, +}\\
\displaystyle{ + \int_{\mathbb{R}} \int_{\partial^*\{u>t\}} (\theta_{\Phi(t)}(x)-1)f(x, n_t(x)^\perp) \;{d} \mathcal{H}^1 {d}t } 
\displaystyle{\,\, + \int_{\mathbb{R}} \int_{\Phi(t)\setminus \partial^*\{u>t\}} \theta_{\Phi(t)}(x)f(x,  n_t(x)^\perp ) \;{d} \mathcal{H}^1 {d} t}.
\end{array}
\end{equation}
The coarea formula in $\BV$ yields that
\begin{equation}\label{Du1}
\int_{\mathbb{R}} \int_{\partial^*\{u>t\}}f(x, n_t(x)^\perp) \;{d}  \mathcal{H}^1 {d}t = \int_{\mathbb{R}^2}f(x, Du(x)^\perp ) \;{d} |Du|. 
\end{equation}
\par The decomposition theorem for the derivative of $\SBV$ functions implies
$$|Du| = |\nabla u| \mathcal{L}^2 +  |D^su| $$
where $\nabla u$ is the approximate gradient of $u$ and $|D^su| =|u^+ - u^-|\mathcal{H}^1\res J_u$ where $J_u$ is the set of approximate jump points of $u$, so
\begin{equation}\label{Du2}
\displaystyle{\int_{\mathbb{R}^2} f(x, Du(x)^\perp  ) \;{d} |Du|= \int_{\mathbb{R}^2} |\nabla u|(x) f\left( x, \nabla u(x)^\perp  \right)  \;{{d} x} + 
\int_{\Omega}  f\left( x,  D^su(x)^\perp \right)  \;{d} |D^su|.}
\end{equation}
\par Moreover  we can consider  the measure $m$ defined as:
\begin{equation}\label{m}
m(A) = \int_{\mathbb{R}} \int_{[\Phi(t)\setminus \partial^*\{u>t\}]\cap A} \theta_{\Phi(t)}(x)\;{d} \mathcal{H}^1\;{d} t + \int_{\mathbb{R}} \int_{ \partial^*\{u>t\}\cap A} (\theta_{\Phi(t)}(x)-1)\;{d} \mathcal{H}^1\;{d} t
\end{equation} 
for every measurable set $A\subset \Omega$ and, from \eqref{ipotbordo}, we get 
\begin{equation}\label{mfront}
 m(\partial \Omega)=0 
\end{equation}
\par Thus, by \eqref{primadecomposizione}, \eqref{Du1}, \eqref{Du2} and \eqref{m}, we can write 
\begin{equation}\label{decomposizionetotale}
\begin{array}{ll}
 \displaystyle{\int_{G_{1}(\Omega)} f(x, S) \;{d} V(x,S) = }& \displaystyle{ \int_{\Omega} |\nabla u|(x) f\left( x, \nabla u(x)^\perp  \right)  \;{{d} x} 
  +\int_{\Omega}  f\left( x,  D^su(x)^\perp \right)  \;{d} |D^su|} \\
 & + \displaystyle{\int_{\Omega } f(x,n_t(x)^\perp) \;{d} m(x) }.
\end{array}
\end{equation}
We deduce from the previous equality that $V$ is a Young varifold.
 In fact, by \eqref{decomposizionetotale}, we get
\begin{equation}\label{YV1}
 \int_{G_{1}(\Omega)} f(x, S) \;{d} V(x,S) = \langle \langle \nu,  |z|f(x,  z^\perp) \rangle\rangle
\end{equation}
where $\nu=(\nu_x, \nu_x^{\infty}, \lambda_\nu)$ is the Young measure defined as:
\begin{equation}\label{YV2}
 \nu_x=\delta_{\nabla u(x)}, \quad \quad  \lambda_\nu =|D^s u| + m, \quad \quad \nu_x^{\infty}=\left\{
\begin{array}{ll}
 \delta_{\frac{D^s u}{|D^s u|}(x)} &  \text{ if\;} x\in \,J_u\\
 \dfrac{1}{2}\left( \delta_{n_t(x)}+ \delta_{- n_t(x)}\right) &  \text{ if\;} x\in\, \text{spt}\;m
\end{array} \right.
\end{equation}
where $n_t(x)$ is a unit normal vector to the system $\Phi(t)$ at $x$. Remark that, as different curves belonging either to the same system or to different systems may intersect only tangentially, the previous measure is well defined.
\par By \eqref{ipotbordo} and \eqref{mfront} we have $\lambda_\nu(\partial \Omega)=0$ and in addition
$$ \int_{\Omega}\int_{\rn}|z| \;{d} \nu_x(z)\;{{d} x}+\lambda_{\nu}(\overline{\Omega}) = \int_{\mathbb{R}}\int_{\Phi(t)} \theta_{\Phi(t)} \;{d} \mathcal{H}^1 \;{d} t\leq G(\Phi)< \infty $$
so, using Theorem \ref{caract}, we get $\nu\in {\bf GY}(\Omega, \mathbb{R}^2)$. In addition it is easy to check that ${\operatorname{Bar}}_\nu =  Du$ so $\nu \in\mbox{{\bf GY}}(u)$.  To end the proof we have to show that $V=V_\nu \in \mathbb{V}(u)$ and 
$$W(V) = \int_{\mathbb{R}} W(\Phi(t)) \;{d} t = G(\Phi).$$
\par For every measurable set $A$,
$$
 \mu_V(A)  = \int_{\mathbb{R}} \int_{\Phi(t)\cap A} \theta_{\Phi(t)} \;{d} \mathcal{H}^1\;{d} t
$$
and for every $X\in \cont^1_0( \Omega, \mathbb{R}^2)$,  by \eqref{azionemiovar} and \eqref{azionelivelli}, we get 
$$
\begin{array}{ll}
\delta V(X) & \displaystyle= \int_{\mathbb{R}} \int_{\Phi(t)} \theta_{\Phi(t)}{\rm{div}}_{\Phi(t)} X \;{d} \mathcal{H}^1\;{d} t=- \int_{\mathbb{R}} \int_{\Phi(t)}  \theta_{\Phi(t)}\langle X, {\bf H}_{V_t} \rangle \;{d} \mathcal{H}^1\;{d} t=- \int_{\Omega} \langle X, {\bf H} \rangle \;{d} \mu_V
\end{array}
$$
where the last equality follows from the locality of the mean curvature (see Theorem \ref{localH}), that guarantees, together with the properties of $\Phi$, that the mean curvature is uniquely defined $\mu_V$-almost everywhere, . 
\par It follows from the expression of the first variation that $ {\bf H}_V = {\bf H}$, thus $\|\delta V\|<\!<\mu_V$ so $V=V_\nu \in \mathbb{V}(u)$. Moreover, by  the coarea formula, 
$$
\begin{array}{ll}
G(\Phi) &= \displaystyle{\int_{\mathbb{R}} W(\Phi(t)) \;{d} t = \int_{\mathbb{R}} \int_{\Phi(t)} \theta_{\Phi(t)} \;{d} \mathcal{H}^1\;{d} t + \int_{\mathbb{R}} \int_{\Phi(t)}  \theta_{\Phi(t)}\arrowvert{\bf H}_{V_t}\arrowvert^p \;{d} \mathcal{H}^1\;{d} t}\\
& \displaystyle{=\mu_V(\Omega) + \int_{\Omega} \arrowvert{\bf H}_{V}\arrowvert^p\;{d} {\mu_V} = W(V).}
\end{array}
$$
\qed

\boss Does a similar result hold in higher dimension ? This is also an open problem. Generalizing to higher dimension the strategy used in dimension $2$ is very delicate because we do not have any description by parametric foliation of our limit functions, see the discussion in the introduction. We believe instead that the Young varifold structure is rich enough and keeps track of sufficiently many information in the limit to allow a desingularization procedure that is necessary to get the representation of the limit Willmore energy using Young varifolds. This is the purpose of ongoing research and we conclude this discussion with a conjecture. 
\eoss
\begin{conj}\label{congettura}
  For every $u\in \BV(\Omega)$ with $\oF(u)<\infty$, there exists $V\in \mathbb{V}(u)$ such that 
$$\oF(u)=W(V).$$
\end{conj}

The set ${\bf GY}(0)$ plays an important role in~\cite{KR}, and it is also of interest in our context. 
\begin{defi}
 Given $u\in \BV(\Omega)$, $\Omega\subset\R^n$, we denote by $\nu_{Du}+ {\bf GY}(0)$ the class of gradient Young measures that decomposes into $\nu_{Du}+\tilde\nu$ for some $\tilde\nu\in {\bf GY}(0)$.
\end{defi}
It is easy to check that $\nu_{Du}+ {\bf GY}(0) \subseteq {\bf GY}(u).$ Remark that if $\nu \in \nu_{Du}+ {\bf GY}(0)$ then 
$$V_\nu = V_{\nu_{Du}}+ V_{\tilde{\nu}}\,, \quad \mbox{with} \quad\tilde{\nu} \in {\bf GY}(0).$$

\begin{prop}\label{struttura}
 The Young measure $\nu$ defined in \eqref{YV2} belongs to $v_{Du}+ {\bf GY}(0)$. 
\end{prop}

\dimo For every $f\in {\bf E}(\Omega; \mathbb{R}^2)$ we have
$$\langle\langle \nu, f\rangle\rangle = \int_{\Omega}f\left( x, \nabla u(x)\right) \;{{d} x} + \int_{\Omega}f^{\infty}\left( x, \dfrac{D^s u}{|D^s u|}\right) \;{d} |D^s u| + \int_{\Omega}f^{\infty}\left( x, n(x)\right) \;{d} m=$$
$$= \langle\langle \nu_{Du}, f\rangle\rangle + \langle\langle \tilde{\nu}, f\rangle\rangle $$
where $\tilde{\nu}$ is defined by the following triplet:
\begin{equation}\label{tildenu}
 \tilde{\nu}_x=\delta_{0}, \quad \quad  \lambda_{\tilde{\nu}} = m, \quad \quad \tilde{\nu}_x^{\infty}= \dfrac{1}{2}\left( \delta_{n(x)}+ \delta_{- n(x)}\right).
\end{equation}
Then $\nu = \nu_{Du}+ \tilde{\nu}$ and it is sufficient to prove that $\tilde{\nu} \in {\bf GY}(0)$. Now,  
 $m(\Omega)< \infty$ and $m(\partial\Omega)= 0$ so, by Theorem \ref{caract}, $\tilde{\nu} \in {\bf GY}(\Omega, \mathbb{R}^2)$. Moreover ${\operatorname{Bar}}_{\tilde{\nu}} = 0$ thus  $\tilde{\nu} \in {\bf GY}(0)$.
\qed

\begin{defi}
 We denote by $\mathbb{V}_0(u)$ the class of Young varifolds $V_\nu\in\mathbb{V}(u)$ such that $\nu \in \nu_{Du}+{\bf GY}(0)$, thus for every $V\in \mathbb{V}_0(u)$
$$V = V_{\nu_{Du}}+ V_{\tilde{\nu}}\,, \quad \mbox{with} \quad\tilde{\nu} \in {\bf GY}(0).$$
\end{defi}

\boss
The definition of $\mathbb{V}_0(u)$ is motivated by the relationship between a general $\nu\in {\bf GY}(u)$ and $\nu_{Du}$. Remark~\ref{reg} shows 
that $\nu_{Du}$ represents $\oF(u)$ for every $u\in \cont^2(\Omega)$ and we have 
$$\oF(u)=W(\nu_{Du}).$$
However, Remark~\ref{non} implies that we cannot expect such a relation for every $u\in \BV(\Omega)$ because  in general $V_{\nu_{Du}}\notin \mathbb{V}(u)$. Thus, a natural question is the following: given $u\in \BV(\Omega)$, how different are the measures $\nu\in {\bf GY}(u)$ and $\nu_{Du}$? 

This leads to characterizing the measures belonging to $\mathbb{V}_0(u)$ in order to estimate $\tilde{\nu}$. Several questions arise naturally. Do $\mathbb{V}_0(u)$ and $\mathbb{V}(u)$ coincide? Does the solution of the minimum problem in 
Proposition \ref{minVu} live in $\mathbb{V}_0(u)$? This would mean that, in order to reach the minimum of $W$, one should take an ``economic`` $\tilde{\nu}$, which is very expectable. But, so far, all these questions remain open.
\eoss

\par Let us examine what can be said about the minimization of $W$ in $\mathbb{V}_0(u)$.

\begin{thm}\label{minYdec}
Let $u\in \BV(\Omega)$, $\Omega\subset\R^n$, such that $\overline{F}(u,\Omega) < \infty$. If there exists $V\in\mathbb{V}_0(u)$ such that $W(v)<\infty$ then the problem 
$$ \Min\; \{W(V):V\in\mathbb{V}_0(u)\}$$
has at least one solution. 
\end{thm}

\dimo Let $\{V_h\} \subset \mathbb{V}_0(u)$ be a minimizing sequence. By definition of $\mathbb{V}_0(u)$ we can take a sequence $\{\tilde{\nu}_h\} \subset {\bf GY}(0)$ such that
$$ V_h = V_{\nu_{Du}}+ V_{\tilde{\nu}_h} \quad \forall h.$$
We can suppose $W(V_h)$ uniformly bounded so we get 
$$\underset{h}{\sup}\; \mu_{V_{\tilde{\nu}_h}}(\Omega)< \infty.$$
Then, by Theorems \ref{comp} and \ref{chiusura2}, there exist a subsequence (not relabeled) of $\{\tilde{\nu}_h\}$ and a Young measure $\tilde{\nu}\in {\bf GY}(0)$ such that
$$ \tilde{\nu}_h \overset{{\bf Y}}{\longrightarrow}  \tilde{\nu}.$$
Then we get $\nu_h=\nu_{Du}+ \tilde{\nu}_h \overset{{\bf Y}}{\longrightarrow} \nu=\nu_{Du}+ \tilde{\nu}$ and $V_h \overset{*}{\rightharpoonup} V_\nu $, where $V_\nu=V_{\nu_{Du}}+ V_{\tilde{\nu}}$ is the Young varifold associated with the gradient Young measure $\nu=\nu_{Du}+ \tilde{\nu}$, $\tilde{\nu}\in {\bf GY}(0)$. 
\par In addition, by Proposition~\ref{236}, we have $\| \delta V_{\nu} \| <\!< \mu_{V_\nu}$ and $W(V_\nu)\leq \underset{n\rightarrow\infty}{\liminf}\; W(V_h).$ Then $V_\nu\in \mathbb{V}_0(u)$ and $W(V_\nu)=\mbox{Min}\;\{W(V):V\in\mathbb{V}_0(u)\}$, and the theorem ensues.
\qed

The next corollary, that follows from Theorem~\ref{sistemaYoung}, points out the relationship in dimension $2$ between Young varifolds in  $\mathbb{V}_0(u)$ and the relaxation problem on a bounded domain. 

\begin{coro}\label{phiyoung}
 Let $u\in \rm{SBV}(\mathbb{R}^2)$ with compact support and such that $\overline{F}(u,\mathbb{R}^2) < \infty$. There exists a bounded open domain $\Omega$ 
 and a Young measure $\nu\in {\bf GY}(\Omega, \mathbb{R}^2)$ such that  $V_\nu \in \mathbb{V}_0(u) $ and 
$$\overline{F}(u,\mathbb{R}^2)=\inf\left\lbrace
\underset{h\rightarrow\infty}{\liminf}\;F(u_h, \Omega) : \{u_h\}\in \cont^2_{\mathrm c}(\Omega), \;
u_h\overset{\lp{1}(\Omega)}{\longrightarrow}u\right\rbrace=W(V_\nu).$$
\end{coro}

\dimo By Theorem \ref{principale} there exists $\Phi\in \mathscr{A}(u)$ such that $\overline{F}(u,\mathbb{R}^2)= G(\Phi)$ and, by Proposition \ref{legameY}, there exists an open bounded domain $\Omega$ with
$$ \bigcup_{t\in \mathbb{R}}(\Phi(t))\subset\subset \Omega$$
and such that $u\in \BV(\Omega)$ and 
$$\inf\left\lbrace
\underset{h\rightarrow\infty}{\liminf}\;F(u_h, \Omega) : \{u_h\}\in \cont^2_{\mathrm c}(\Omega), \;
u_h\overset{\lp{1}(\Omega)}{\longrightarrow}u\right\rbrace=\overline{F}(u,\mathbb{R}^2)=G(\Phi)$$
\par Then, using Theorem~\ref{sistemaYoung}, we can define a suitable Young measure $\nu\in {\bf GY}(\Omega, \mathbb{R}^2)$ such that $V_\nu\in \mathbb{V}(u)$ and 
$$\inf\left\lbrace
\underset{h\rightarrow\infty}{\liminf}\;F(u_h, \Omega) : \{u_h\}\in \cont^2_{\mathrm c}(\Omega), \;
u_h\overset{\lp{1}(\Omega)}{\longrightarrow}u\right\rbrace=G(\Phi)=W(V_\nu).$$
Moreover, by  Proposition \ref{struttura}, $V_\nu \in \mathbb{V}_0(u)$ and, by Theorem  \ref{minYdec}, we have
$$\Min \{W(V):V\in\mathbb{V}_0(u)\} \leq \inf\left\lbrace
\underset{h\rightarrow\infty}{\liminf}\;F(u_h, \Omega) : \{u_h\}\in \cont^2_{\mathrm c}(\Omega), \;
u_h\overset{\lp{1}(\Omega)}{\longrightarrow}u\right\rbrace.$$
\qed
\boss
In the previous corollary, the quantity 
$$\inf\left\lbrace
\underset{h\rightarrow\infty}{\liminf}\;F(u_h, \Omega) : \{u_h\}\in \cont^2_{\mathrm c}(\Omega), \;
u_h\overset{\lp{1}(\Omega)}{\longrightarrow}u\right\rbrace$$
corresponds to the relaxation of $F$ using approximating functions in $\cont^2_{\mathrm c}$ instead of $\cont^2$. This definition is well
 posed but induces different properties for the relaxed functional, as discussed in~\cite{MN} where some examples are also provided. We shall adopt this definition over this remark to show the link between Young varifolds, Young measures and the Willmore functional.
   
\par In the regular planar case (i.e. $u\in \cont^2_{\mathrm c}(\Omega)$, $\Omega\subset\R^2$) there are two natural frameworks to represent $\overline{F}$:
 \begin{itemize}
  \item by a coarea-type formula, using the representation by systems of curves of class $\wpq{2}{p}$:
 $$ \overline{F}(u, \Omega) = G(\Phi[u]) \quad\quad \forall u\in \cont^2_{\mathrm c}(\Omega);$$
 \item by the varifold theory (see Remark \ref{reg}), using Young varifolds:
 $$ \overline{F}(u, \Omega) = W(V_{\nu_{Du}})\quad\quad \forall u\in \cont^2_{\mathrm c}(\Omega).$$
 \end{itemize}
\par Theorem~\ref{sistemaYoung} and Corollary~\ref{phiyoung} state that, as in the regular case, Young varifolds provide a natural framework to represent $\overline{F}$,  at least in dimension $2$. Moreover,  the definition of the Young measure $\nu$ in \eqref{YV2} shows the relationship between the Young varifold representing $\overline{F}(u,\Omega)$ and $V_{\nu_{Du}}$: the Young varifold involves an additional term that contains in particular all "ghost" parts, as in Figure~\ref{cusp3}.
\eoss

\section{Conclusion}
 We introduced in this paper a new framework to address the relaxation of a generalized Willmore functional. We believe that this combination of Young measures and varifolds is the right approach to track, in the limit of oscillations and concentration, the behavior of the energy. In addition, this framework has a major advantage over representations by foliation: the compactness and the semicontinuity of the energy (under some constraints) come easily, as shown in Corollary~\ref{corocomp}
 \par 
We showed in the paper several properties of Young varifolds, we proved a representation result for $\oF$ in dimension 2 (Theorem~\ref{sistemaYoung}, Corollary~\ref{phiyoung}), and, in any dimension $\geq 2$, we proved an inequality that involves a minimum problem for  Young varifolds with prescribed barycenter (Theorem~\ref{minVu}).
 
 There are several obstacles to get a full understanding of the problem:
 \begin{itemize}
 \item The class of Young varifolds associated with a given function is very rich. This is due to the fact that there are infinitely many Young measures that are the limits of sequences of gradient Young measures associated with smooth functions $u_h$ that converge weakly-* to $0$ in $\BV$. It is reasonable to think, however, that minimizing $W$ in $\mathbb{V}(u)$ reduces considerably the measures of interest.
 \item The proof of the representation result in dimension 2 (Theorem~\ref{sistemaYoung}, Corollary~\ref{phiyoung}) strongly relies on Theorem~\ref{principale} that involves a curve stretching technique. This technique can hardly be generalized to higher dimensions, in particular because the accurate description of the boundaries of sets with finite relaxed energy is still an open problem. Another strategy, that has been totally unexplored so far, requires understanding how a limit Young varifold can be regularized with a control of the energy, using in particular the directions of concentration indicated by $\nu_x^\infty$.
 \end{itemize}
 Beyond the relaxation of the generalized Willmore functional, we may think at other problems that could be tackled with Young varifolds, for instance understanding precisely $\Gamma$-limits of suitable functionals when the underlying function is unsmooth, or defining accurately the flow associated with the generalized Willmore functional in order to have a new look at the critical points. We believe that the versatility of Young varifolds makes them delicate but powerful tools.
 
 \medskip
\subsection*{Acknowledgments}
We warmly thank Vicent Caselles and Matteo Novaga for the many discussions we had on the relaxation of the Willmore functional before the current work was initiated. We also thank Giovanni Bellettini and Luca Mugnai for our discussions on related phase-field approximation issues.\\
This work was supported by the French "Agence Nationale de la Recherche" (ANR), under grant FREEDOM (ANR07-JCJC-0048-01).

\bibliographystyle{plain}
\bibliography{biblio-new}

\end{document}